\documentclass{tran-l}
 \newtheorem{thm}{Theorem}[subsection]

 \newtheorem{prop}[thm]{Proposition}
 \theoremstyle{definition}
 \newtheorem{defn}[thm]{Definition}
 \theoremstyle{remark}
 \newtheorem{rem}[thm]{Remark}
 \numberwithin{equation}{subsection}
 \newtheorem{exa}[thm]{Example}

\usepackage{graphicx}
\begin{document}

\title[Chen-Ruan Cohomology ]
 {The Chen-Ruan Cohomology of Weighted Projective
Spaces }

\author{Yunfeng Jiang }

\address{Institute of Systems Science, Academy of Mathematics and System
Sciences, Chinese Academy of Sciences, Beijing, 100080, P.R.China}

\email{jiangyf@mail.amss.ac.cn}

\thanks{}

\thanks{}



\date{}

\dedicatory{}

\commby{}


\begin{abstract}
Chen and Ruan [6] defined a very interesting cohomology theory for
orbifolds, which is now called Chen-Ruan cohomology. The primary
objective of this paper is to compute the Chen-Ruan cohomology
rings of the weighted projective spaces,  a class of important
spaces in  physics. The classical tools (Chen-Ruan cohomology,
toric varieties, the localization technique) which have been
proved to be successful are used to study the orbifold cohomology
of weighted projective spaces. Given a weighted projective space
${\bf P}^{n}_{q_{0}, \cdots, q_{n}}$, we determine all of its
twisted sectors and the corresponding degree shifting numbers, and
we calculate the orbifold cohomology group of ${\bf P}^{n}_{q_{0},
\cdots, q_{n}}$. For a general reduced weighted projective space,
we give a formula to compute the 3-point function which is the key
in the definition of Chen-Ruan cohomology ring. Finally we
concretely calculate the Chen-Ruan cohomology ring of  weighted
projective space ${\bf P}^{5}_{1,2,2,3,3,3}$.
\end{abstract}

\maketitle

\section{Introduction}
The notion of Chen-Ruan orbifold cohomology has appeared in
physics as a result of studying the string theory on global
quotient orbifold, (see [9] and [10]). In addition to the usual
cohomology of the global quotient, this space included the
cohomology of so-called twisted sectors. Zaslow [25] gave a lot of
examples of global quotients and computed their orbifold
cohomology spaces.  But the real mathematical definition of
orbifold cohomology was given by Chen and Ruan [6] for arbitrary
orbifolds. The most interesting feature of this new cohomology
theory, besides the generalization of non global quotients, is the
existence of  a ring structure which was previously missing. This
ring structure is obtained from Chen-Ruan's orbifold quantum
cohomology construction (see [7]) by restricting to the class
called ghost maps, the same as the ordinary cup product may be
obtained by quantum cup product. Since the Chen-Ruan cohomology
appeared, the problem of how to calculate the orbifold cohomology
has been considered by several authors. Chen and Ruan [6] gave
several simple examples. Chen Hao [4] computed the orbifold
cohomology group of moduli space $\mathcal{M}$$_{0,n}$$/S_{n}$.
B.Doug Park and Mainak Poddar [22] considered the Chen-Ruan
cohomology ring of the mirror quintic. All the above examples are
orbifold global quotients. In this paper we calculate the
Chen-Ruan  cohomology rings of weighted projective spaces-a large
class of non-global quotient orbifolds.

A very power tool to compute the Chen-Ruan cohomology of weighted
projective spaces is the method of toric varieties. The theory of
toric varieties establishes a classical connection between
algebraic geometry and the theory of convex polytopes. From the
fan of a toric variety, we can obtain a lot of information about
the toric variety. In particular, when the fan $\Sigma$ of a toric
variety $X$ is simplicial, the toric variety $X$ is an orbifold
with finite abelian groups as local groups. In this paper, we take
the weighted projective spaces as simplicial toric varieties with
local isotropy groups the finite cyclic groups. And then using the
properties of toric varieties induced from the fans, we calculate
the Chen-Ruan  cohomology group of any weighted projective space.

To calculate the Chen-Ruan cohomology ring of the weighted
projective space, we use the Riemann bilinear relations for
periods [15] to identify the obstruction bundle. Up to now, except
that the obstruction bundle for the mirror quintic example
calculated by Park and Poddar [22] is a nontrivial line bundle,
all the other calculated obstruction bundles of the examples are
trivial. In this paper, the obstruction bundle we consider is the
Whitney sum of some line bundles, generalizing the case of mirror
quintic example. And we also introduce the localization techniques
[3] which should work for toric varieties to compute the 3-point
function which is the key in the orbifold cup product [6]. In
particular, we give a concrete example.

On the other hand, a very interesting aspect of calculating the
Chen-Ruan cohomology rings of weighted projective spaces lies in a
conjecture of Ruan. In string theory, physicists suggest that the
orbifold string theory of an orbifold should be equivalent to the
ordinary string theory of its crepant resolution. For the orbifold
cohomology, the Cohomology Hyperkahler Resolution Conjecture of
Ruan (see [23]) states that the Chen-Ruan cohomology ring of an
orbifold should be isomorphic to the ordinary cohomology ring of
its hyperkahler resolution. I hope that my calculation of the
Chen-Ruan cohomology ring of the weighted projective space may
contribute to this interesting problem.

The thesis is outlined as follows. Section 2 is a review of some
basic facts concerning orbifold, Chen-Ruan cohomology and
simplicial toric varieties.  In section 3 we introduced the basic
concept of the weighted projective space. In section 4 we discuss
the Chen-Ruan cohomology group of any weighted projective space.
And in the section 5 we compute the ring structure of the
Chen-Ruan cohomology of the weighted projective space.


\section{Preliminaries}
\subsection{Orbifold and orbifold vector bundle.}
\begin{defn}An orbifold structure on a
Hausdorff,  separate topological space $X$ is given by an open
cover $\mathcal{U}$ of $X$ satisfying the following conditions.

(1) Each element $U$ in $\mathcal{U}$ is uniformized, say by
$(V,G,\pi)$.  Namely, $V$ is a smooth manifold and $G$ is a finite
group acting smoothly on $V$ such that $U=V/G$ with $\pi$ as the
quotient map. Let $Ker(G)$ be the subgroup of $G$ acting trivially
on $V$.

(2) For $U^{'}\subset U$,  there is a collection of injections
$(V^{'},G^{'},\pi^{'})\longrightarrow (V,G,\pi)$. Namely, the
inclusion $i:U^{'}\subset U$ can be lifted to maps $\widetilde{i}:
V^{'}\longrightarrow V $ and an injective homomorphism $i_{*}:
G^{'}\longrightarrow G$ such that $i_{*}$ is an isomorphism from
$Ker(G^{'})$ to $Ker(G)$ and  $\widetilde{i}$ is
$i_{*}$-equivariant.

(3) For any point $x\in U_{1}\cap U_{2}$, $U_{1},U_{2}\in
\mathcal{U} $, there is a $U_{3}\in \mathcal{U} $ such that $x\in
 U_{3} \subset U_{1}\cap U_{2}$.
\end{defn}

For any point $x\in X$, suppose that $(V,G,\pi)$ is a uniformizing
neighborhood and $\overline{x}\in \pi^{-1}(x)$. Let $G_{x}$ be the
stabilizer of $G$ at $\overline{x}$. Up to conjugation, it is
independent of the choice of $\overline{x}$ and is called the
$local~ group$ of $x$. Then there exists a sufficiently small
neighborhood $V_{x}$ of $\overline{x}$ such that
$(V_{x},G_{x},\pi_{x})$ uniformizes a small neighborhood of $x$,
where $\pi_{x}$ is the restriction $\pi\mid V_{x}$.
$(V_{x},G_{x},\pi_{x})$ is called $a$ $local$ $ chart$ at $x$. The
orbifold structure is called $reduced $ if the action of $G_{x}$
is effective for every $x$.

Let $pr: E\longrightarrow X$ be a rank $k$ complex $orbifold~
bundle$ over an orbifold $X$([6]). Then a uniformizing  system for
$E\mid U=pr^{-1}(U)$ over a uniformized subset $U$ of $X$ consists
of the following data:

(1) A uniformizing system $(V,G,\pi)$ of $U$.

(2) A uniformizing system $(V \times \bf{C^{k}},$ $ G,
\widetilde{\pi})$ for $E\mid U$. The action of $G$ on $V\times
\mathbf{C}^{k}$ is an extension of the action of $G$ on $V$ given
by $g\cdot(x,v)=(g\cdot x,\rho(x,g)v)$ where $\rho: V\times
G\longrightarrow Aut(\mathbf{C}^{k})$ is a smooth map satisfying:
$$\rho(g\cdot x,h)\circ \rho(x,g)=\rho(x,hg),g,h\in G, x\in V.$$

(3)The natural projection map $\widetilde{pr}: V \times
\bf{C}^{k}\longrightarrow V$ satisfies $\pi\circ
\widetilde{pr}=pr\circ \widetilde{\pi}$.

By an orbifold connection $\bigtriangleup$ on $E$ we mean an
equivariant connection that satisfies
$\bigtriangleup=g^{-1}\bigtriangleup g$ for every uniformizing
system of $E$. Such a connection can be always obtained by
averaging an equivariant partition of unity.

\subsection{Twisted Sectors and Chen-Ruan  Cohomology.}

The most physical idea is twisted sectors. Let $X$ be an orbifold.
Consider the set of pairs:

$$\widetilde{X}_{k}=\{(p,(\mbox{\bf{g}})_{G_{p}})|p\in X,
\mbox{\bf{g}}=(g_{1},\cdots,g_{k}),g_{i}\in G_{p}\}$$ where
$(\mbox{\bf{g}})_{G_{p}}$ is the conjugacy class of $k$-tuple
$\mbox{\bf{g}}=(g_{1},\cdots,g_{k})$ in $G_{p}$. We use $G^{k}$ to
denote the set of $k$-tuples. If there is no confusion, we will
omit the subscript $G_{p}$ to simplify the notation. Suppose that
$X$ has an orbifold structure $\mathcal{U}$ with uniformizing
systems $(\widetilde{U},G_{U},\pi_{U})$. From Chen and Ruan [6],
also see [18],  we have: $\widetilde{X}_{k}$ is naturally an
orbifold,  with the generalized orbifold structure at
$(p,(\mbox{\bf{g}})_{G_{p}})$ given by
$(V_{p}^{\bf{g}},C(\mbox{\bf{g}}),\pi:V_{p}^{\bf{g}}\longrightarrow
V_{p}^{\bf{g}}/ C(\mbox{\bf{g}}))$,  where
$V_{p}^{\bf{g}}=V_{p}^{g_{1}}\cap \cdots V_{p}^{g_{k}}$,
$C(\mbox{\bf{g}})=C(g_{1})\cap\cdots C(g_{k})$. Here
$\mbox{\bf{g}}=(g_{1}, \cdots, g_{k})$,  $V_{p}^{g}$ stands for
the fixed point set of $g$ in $V_{p}$. When $X$ is almost complex,
$\widetilde{X_{k}}$ inherits an almost complex structure from $X$,
and when $X$ is closed,  $\widetilde{X_{k}}$ is finite disjoint
union of closed orbifolds.

Now we describe the the connected components of
$\widetilde{X_{k}}$,  Recall that every point $p$ has a local
chart $(V_{p},G_{p},\pi_{p})$ which gives a local uniformized
neighborhood $U_{p}=\pi_{p}(V_{p})$. If $q\in U_{p}$,  up to
conjugation there is a unique injective homomorphism $i_{*}:
G_{q}\longrightarrow G_{p}$.  For $\mbox{\bf{g}}\in (G_{q})^{k}$,
the conjugation class $i_{*}(\mbox{\bf{g}})_{q}$ is well defined.
We define an equivalence relation $i_{*}(\mbox{\bf{g}})_{q}\cong
(\mbox{\bf{g}})_{q}$. Let $T_{k}$ denote the set of equivalence
classes.To abuse the notation,  we use $(\mbox{\bf{g}})$ to denote
the equivalence class which $(\mbox{\bf{g}})_{q}$ belongs to.  We
will usually denote an element of $T_{1}$ by $(g)$. It is clear
that $\widetilde{X_{k}}$ can be decomposed as a disjoint union of
connected components:

$$\widetilde{X_{k}}=\bigsqcup_{({\bf g}) \in T_{k}}X_{(\bf{g})}$$

Where $X_{(\bf{g})}=\{(p,(\mbox{\bf{g}}^{'})_{p})|
\mbox{\bf{g}}^{'} \in (G_{p})^{k},(\mbox{\bf{g}}^{'})_{p}\in
(\mbox{\bf{g}})\}$. Note that for ${\bf g}=(1,\cdots,1)$, we have
$X_{(\bf{g})}=X$.  A component $X_{(\bf{g})}$ is called a
$k-multisector$, if $\bf{g}$ is not the identity. A component of
$X_{(g)}$ is simply called a $twisted$ $sector$. If $X$ has an
almost complex, complex or kahler structure,  then $X_{(\bf{g})}$
has the analogous structure induced from $X$. We define
$$T_{3}^{0}=\left\{({\bf g})=(g_{1},g_{2},g_{3})\in T_{3}| g_{1}g_{2}g_{3}=1\right\}.$$
Note that there is an  one to one correspondence between $T_{2}$
and $T_{3}^{0}$ given by $(g_{1},g_{2})\longmapsto
(g_{1},g_{2},(g_{1}g_{2})^{-1})$.

Now we define the Chen-Ruan cohomology. Assume that $X$ is a
$n$-dimensional compact almost complex orbifold with almost
structure $J$. Then for a point $p$ with nontrivial group $G_{p}$,
$J$ gives rise to an effective representation
$\rho_{p}:G_{p}\longrightarrow GL(n,\bf{C})$. For any $g\in
G_{p}$, we write $\rho_{p}(g)$, up to conjugation, as a diagonal
matrix
$$diag\left(e^{2\pi i\frac{m_{1,g}}{m_{g}}},\ldots,e^{2\pi i\frac{m_{n,g}}{m_{g}}}\right).$$
where $m_{g}$ is the order of $g$ in $G_{p}$, and $0\leq m_{i,g}<
m_{g}$. Define a function $\iota: \widetilde{X_{1}}\longrightarrow
\bf{Q} $ by
$$
\iota(p,(g)_{p})=\sum_{i=1}^{n}\frac{m_{i,g}}{m_{g}}.
$$
We can see that the function $\iota:
\widetilde{X_{1}}\longrightarrow \bf{Q} $ is locally constant and
$\iota=0$ if $g=1$. Denote  its value on $X_{(g)}$ by $\iota_{g}$.
We call $\iota_{g}$ the degree shifting number of $X_{(g)}$. It
has the following properties:

(1) $\iota_{(g)}$ is an integer iff $\rho_{p}(g)\in SL(n,\bf{C})$;

(2)
$\iota_{(g)}+\iota_{(g^{-1})}=rank(\rho_{p}(g)-Id)=n-dim_{\bf{C}}X_{(g)}$.\\

A $C^{\infty}$ differential form on $X$ is a $G$-invariant
differential form on $V$ for each uniformizing system $(V,G,\pi)$.
Then orbifold integration is defined as follows. Suppose $U=V/G$
is connected, for any compactly supported differential $n$-form
$\omega$ on $U$, which is, by definition, a $G$-invariant $n$-form
$\widetilde{\omega}$ on $V$,

$$
 \int_{U}^{orb}\omega
:=\frac{1}{|G|}\int_{V}\widetilde{\omega}  \eqno{(2.1)}
$$

Where $|G|$ is the order of $G$.  The orbifold integration over
$X$ is defined by using a $C^{\infty}$ partition of unity. The
orbifold  integration coincides with the usual measure theoretic
integration iff the orbifold structure is reduced.

Holomorphic forms for a complex orbifold $X$ are again obtained by
patching $G$-invariant holomorphic forms on the uniformizing
system $(V,G,\pi)$. We consider the Cech cohomology groups of $X$
and $X_{(\bf{g})}$ with coefficients in the sheaves of holomorphic
forms. The Cech cohomology groups can be identified with the
Dolbeault cohomology groups of $(p,q)$-forms [2].

\begin{defn}([6]) Let $X$ be a closed
complex orbifold, we define the orbifold cohomology group of $X$
by

$$H_{orb}^{d}:=\bigoplus_{(g)\in
T_{1}}H^{d-2\iota_{(g)}}(X_{(g)},\bf{Q})$$

For $0\leq p,q\leq dim_{\bf{C}}X$, we define the orbifold
Dolbeault cohomology group of $X$ by

$$H_{orb}^{p,q}(X):=\bigoplus_{(g)\in
T_{1}}H^{p-\iota_{(g)},q-\iota_{(g)}}(X_{(g)},\bf{C})$$
\end{defn}

\subsection{The Obstruction Bundle.}
Choose $(\bf{g})$$=(g_{1},g_{2},g_{3})\in T_{3}^{0}$. Let
$(p,({\bf g})_{p})$ be a generic point in $X_{(\bf{g})}$. Let
$K({\bf g})$ be the subgroup of $G_{p}$ generated by $g_{1}$ and
$g_{2}$. Consider an orbifold Riemann sphere with three orbifold
points $(S^{2},(p_{1},p_{2},p_{3}),(k_{1},k_{2},$ $k_{3}))$. When
there is no confusion, we will simply denote it by $S^{2}$. The
orbifold fundamental group is:

$$\pi_{1}^{orb}(S^{2})=\{\lambda_{1},\lambda_{2},\lambda_{3}|\lambda_{i}^{k_{i}}=1,\lambda_{1}\lambda_{2}\lambda_{3}=1\}$$

Where $\lambda_{i}$ is represented by a loop around the marked
$p_{i}$. There is a surjective homomorphism
$$\rho: \pi_{1}^{orb}(S^{2})\longrightarrow K(\bf{g})$$
specified by mapping $\lambda_{i}\longmapsto g_{i}$. $Ker(\rho)$
is a finite-index subgroup of $\pi_{1}^{orb}(S^{2})$. Let
$\widetilde{\Sigma}$ be the orbifold universal cover of $S^{2}$.
Let $\Sigma=\widetilde{\Sigma}/Ker(\rho)$. Then $\Sigma$ is
smooth,  compact and $\Sigma/K({\bf g})=S^{2}$. The genus of
$\Sigma$ can be computed using Riemann Hurwitz formula for Euler
characteristics of a branched covering,  and turns out to be

$$g(\Sigma)=\frac{1}{2}(2+|K(\mbox{\bf{g}})|-\Sigma_{i=1}^{3}\frac{|K(\bf{g})|}{k_{i}}) \eqno{(2.2)}$$

$K(\bf{g})$ acts holomorphically on $\Sigma$  and hence
$K(\bf{g})$ acts on $H^{0,1}(\Sigma)$. The "obstruction bundle"
$E_{(\bf{g})}$ over $X_{(\bf{g})}$ is constructed as follows. On
the local chart $(V_{p}^{\bf{g}},C(\bf{g}),\pi)$ of
$X_{(\bf{g})}$, $E_{(\bf{g})}$ is given by $(TV_{p}\otimes
H^{0,1}(\Sigma))^{K(\bf{g})}\times V_{p}^{\bf{g}}\longrightarrow
V_{p}^{\bf{g}}$, where $(TV_{p}\otimes
H^{0,1}(\Sigma))^{K(\bf{g})}$ is the $K(\bf{g})$-invariant
subspace. We define an action of  $C(\bf{g})$ on $TV_{p}\otimes
H^{0,1}(\Sigma)$, which is the usual one on $TV_{p}$ and trivial
on $H^{0,1}(\Sigma)$. The the action of $C(\bf{g})$ and
$K(\bf{g})$ commute and $(TV_{p}\otimes
H^{0,1}(\Sigma))^{K(\bf{g})}$ is invariant under $C(\bf{g})$. Thus
we have obtained an action of $C(\bf{g})$ on $(TV_{p}\otimes
H^{0,1}(\Sigma))^{K(\bf{g})}\times V_{p}^{\bf{g}}\longrightarrow
V_{p}^{\bf{g}}$, extending the usual one on $V_{p}^{\bf{g}}$.
These trivializations fit together to define the bundle
$E_{(\bf{g})}$ over $X_{(\bf{g})}$. If we set $e:
X_{(\bf{g})}\longrightarrow X$ to be the map given by $(p,({\bf
g})_{p})\longmapsto p$, one may think of $E_{(\bf{g})}$ as
$(e^{*}TX\otimes H^{0,1}(\Sigma))^{K(\bf{g})}$. The rank of
$E_{(\bf{g})}$  is given by the formula [6]:

$$rank_{\bf{C}}(E_{(\bf{g})})=dim_{\bf{C}}(X_{(\bf{g})})-dim_{\bf{C}}(X)+\Sigma_{j=1}^{3}\iota_{(g_{j})}
\eqno{(2.3)}$$

\subsection{Orbifold cup product.} First, there is a
natural map $I: X_{(g)}\longrightarrow X_{(g^{-1})}$ defined by
$(p,(g)_{p})\longmapsto (p,(g^{-1})_{p})$.

\begin{defn}Let $n=dim_{\bf{C}}(X)$. For any
integer $0\leq n\leq 2n$, the pairing

$$<,>_{orb}: H^{d}_{orb}(X)\times H^{2n-d}_{orb}(X)\longrightarrow \bf{Q}$$
is defined by taking the direct sum of

$$<,>_{orb}^{(g)}: H^{d-2\iota_{(g)}}(X_{(g)};\mbox{\bf{Q}})\times H^{2n-d-2\iota_{(g^{-1})}}(X_{(g^{-1})};\bf{Q})
\longrightarrow \bf{Q}$$ where
$$<\alpha,\beta>_{orb}^{(g)}=\int_{X_{(g)}}^{orb}\alpha\wedge
I^{*}(\beta)$$
\end{defn}

for $\alpha\in H^{d-2\iota_{(g)}}(X_{(g)};\mbox{\bf{Q}})$, and
$\beta\in H^{2n-d-2\iota_{(g^{-1})}}(X_{(g^{-1})};\bf{Q})$.

Choose an orbifold connection $A$ on $E_{(\bf{g})}$. Let
$e_{A}(E_{(\bf{g})})$ be the Euler form computed from the
connection $A$ by Chen-Weil theory. Let $\eta_{j}\in
H^{d_{j}}(X_{(g_{j})};\mbox{\bf{Q}})$, for $j=1,2,3$. Define maps
$e_{j}: X_{(\bf{g})}\longrightarrow X_{(g_{j})}$ by
$(p,(\mbox{\bf{g}})_{p})\longmapsto (p,(g_{j})_{p})$.

\begin{defn} Define the 3-point function to be

$$<\eta_{1},\eta_{2},\eta_{3}>_{orb}:=\int_{X_{(\bf{g})}}^{orb}e^{*}_{1}\eta_{1}\wedge
 e^{*}_{2}\eta_{2}\wedge e^{*}_{3}\eta_{3}\wedge e_{A}(E_{(\bf{g})})
 \eqno{(2.4)}$$
\end{defn}

Note that the above integral does not depend on the choice of $A$.
As in the definition 2.4.1, we extend the 3-point function to
$H_{orb}^{*}(X)$ by linearity. We define the orbifold cup product
by the relation

$$<\eta_{1}\cup_{orb}\eta_{2},\eta_{3}>_{orb}:=<\eta_{1},\eta_{2},\eta_{3}>_{orb}  \eqno{(2.5)}$$

Again we extend $\cup_{orb}$ to $H_{orb}^{*}(X)$ via linearity.
Note that if $({\bf g})=(1,1,1)$, then
$\eta_{1}\cup_{orb}\eta_{2}$ is just the ordinary cup product
$\eta_{1}\cup \eta_{2}$ in $H^{*}(X)$.

\subsection{Simplicial Toric Varieties as Orbifolds.}

A toric variety is a normal variety with an action of an algebraic
torus which admits an open dense orbit homeomorphic  to the torus.
Every toric variety is described by a set of combinatoric data,
called a fan $\Xi$ in a lattice $N$, [13], [20]. $\Xi$ is
simplicial if every cone $\sigma$ in $\Xi$ is generated by a
subset of a basis of ${\bf R}^{n}=N\otimes \bf{R}$.

We now describe  the orbifold structure of simplicial toric
varieties, see Poddar [20]. Let $\Xi$ be any simplicial fan in a
$n$ dimensional lattice $N$. $X_{\Xi}$ be the corresponding toric
variety. For a cone $\tau\in \Xi$, denote the set of its primitive
1-dimensional generators by $\tau[1]$, the corresponding affine
open subset of $X_{\Xi}$ by $U_{\tau}$, and the corresponding
torus orbit by $O_{\tau}$. We write $\nu\leq \tau$ if the cone
$\nu$ is a face of the cone $\tau$, and $\nu< \tau$ if it is a
proper subface. $U_{\tau}=\sqcup_{\nu\leq \tau}O_{\nu}$. Let $M=
Hom(N,\bf{Z})$ be the dual lattice of $N$ with dual pair $<,>$.
For any cone $\tau\in \Xi$, denote its dual cone in $M\otimes
\bf{R}$ by $\check{\tau}$. Let $S_{\tau}=\check{\tau}\cap M$.
$\bf{C}(S_{\tau})$ is the $\bf{C}$-algebra with generators
$\chi^{m}$ for each $m\in S_{\tau}$ and relation
$\chi^{m}\chi^{m^{'}}=\chi^{m+m^{'}}$.
$U_{\tau}=Spec(\mbox{\bf{C}}[S_{\tau}])$.

Then the orbifold structure of the toric variety $X_{\Xi}$ can be
described as follows. Let $\sigma$ be any $n$ dimensional cone of
$\Xi$. Let $v_{1},\cdots,v_{n}$ be the primitive 1 dimension
generstors of $\sigma$. These are linearly independent in
$N_{\bf{R}}=N\otimes \bf{R}$. Let $N_{\sigma}$ be the sublattice
of $N$ generated by $v_{1},\cdots,v_{n}$. And let
$G_{\sigma}=N/N_{\sigma}$ be the quotient group, then $G_{\sigma}$
is finite and abelian.

Let $\sigma^{'}$ be the cone $\sigma$ regarded in $N_{\sigma}$.
Let $\check{\sigma}^{'}$ be the dual cone of  $\sigma^{'}$ in
$M_{\sigma}$, the dual lattice of $N_{\sigma}$.
$U_{\sigma^{'}}=spec(\mbox{\bf{C}}[\check{\sigma}^{'}\cap
M_{\sigma}])$. Note that $\sigma^{'}$ is a smooth cone in
$N_{\sigma}$. So $U_{\sigma^{'}}\cong {\bf C}^{n}$.

Now there is a canonical dual pairing $M_{\sigma}/M\times
N/N_{\sigma}\longrightarrow \bf{Q/Z}\longrightarrow \bf{C^{*}}$,
the first map by the pairing $<,>$ and the second by $q\longmapsto
exp(2\pi iq)$.  Now $G_{\sigma}$ acts on
$\mbox{\bf{C}}[M_{\sigma}]$, the group ring of $M_{\sigma}$, by:
$v(\chi^{u})=exp(2\pi i<u,v>)\chi^{u}$, for $v\in N$ and $u\in
M_{\sigma}$. Note that

$$(\mbox{\bf{C}}[M_{\sigma}])^{G_{\sigma}}=\mbox{\bf{C}}[M]
\eqno{(2.6)}$$

Thus $G_{\sigma}$ acts on $U_{\sigma^{'}}$. Let $\pi_{\sigma}$ be
the quotient map. Then $U_{\sigma}=U_{\sigma^{'}}/G_{\sigma}$. So
$U_{\sigma}$ is uniformized by
$(U_{\sigma^{'}},G_{\sigma},\pi_{\sigma})$. For any $\tau<
\sigma$, the orbifold structure on $U_{\tau}$ is the same as the
one induced from the uniformizing system on $U_{\sigma}$. Then by
the description of the toric gluing it is clear that
$\{(U_{\sigma^{'}},G_{\sigma},\pi_{\sigma}): \sigma\in \Xi[n]\}$
defines a reduced orbifold structure on $X_{\Xi}$. We give a more
explicit verification of this fact below.

Let $B$ be the nonsingular matrix with generators
$v_{1},\cdots,v_{n}$ of $\sigma$ as rows. Then
$\check{\sigma}^{'}$ is generated in $M_{\sigma}$ by the column
vectors $v^{1},\cdots,v^{n}$ of the matrix $B^{-1}$. So
$\chi^{v^{1}},\cdots,\chi^{v^{n}}$ are the coordinates of
$U_{\sigma^{'}}$. For any $k=(k_{1},\cdots,k_{n})\in N$, the
corresponding coset $[k]\in G_{\sigma}$ acts on $U_{\sigma^{'}}$
in these coordinates as a diagonal matrix: $diag(exp(2\pi
ic_{1}),\cdots,exp(2\pi ic_{n}))$, where $c_{i}=<k,v^{i}>$. Such a
matrix is uniquely represented by an $n$-tuple
$a=(a_{1},\cdots,a_{n})$ where $a_{i}\in [0,1)$ and
$c_{i}=a_{i}+b_{i}, b_{i}\in \bf{Z}$. In matrix notation,
$kB^{-1}=a+b \cong k=aB+bB$. We denote the integral vector $aB$ in
$N$ by $k_{a}$ and the diagonal matrix corresponding to $a$ by
$g_{a}$. $k_{a}\longleftrightarrow g_{a}$ gives a one to one
correspondence between the elements of $G_{\sigma}$ and the
integral vector in $N$ that are linear combinations  of the
generators of $\sigma$ with coefficient in $[0,1)$.

Now let us examine the orbifold chart induced by
$(U_{\sigma^{'}},G_{\sigma},\pi_{\sigma})$ at any point $p\in
U_{\sigma}$. By the orbit decomposition, there is a unique
$\tau\in \Xi$ such that $p\in O_{\tau}$. We assume $\tau$ is
generated by $v_{1},\cdots,v_{j},j\leq n$. Then any preimage of
$p$ with respect to $\pi_{\sigma}$ has coordinates
$\chi^{v^{i}}=0$ iff $i\leq j$. Let
$z=(0,\ldots,0,z_{j+1},\ldots,z_{n})$ be one such preimage. Let
$G_{\tau}:=\{g_{a}\in G_{\sigma}: a_{i}=0 ~if ~j+1\leq i\leq n
\}$. We can find a small neighborhood $W\subset ({\bf
C}^{*})^{n-j}$ of $(z_{j+1},\ldots,z_{n})$ such that the
inclusions ${\bf C}^{j}\times W\hookrightarrow U_{\sigma^{'}}$ and
$G_{\tau}\hookrightarrow G_{\sigma}$ induces an injection of
uniformizing systems $({\bf C}^{j}\times W,
G_{\tau},\pi)\hookrightarrow
(U_{\sigma^{'}},G_{\sigma},\pi_{\sigma})$ on some small open
neighborhood $U_{p}$ of $p$. So we have $G_{p}=G_{\tau}$ and an
orbifold chart $({\bf C}^{j}\times W, G_{\tau},\pi)$. Note that
$G_{\tau}$ can be constructed from the set
$\{k_{a}=\Sigma_{i=1}^{j}a_{i}v_{i}: k_{a}\in N, a_{i}\in [0,1)\}$
which is completely determined by $\tau$ and hence is independent
of $\sigma$.


\section{The Weighted Projective Spaces}
\subsection{The Definition and the Orbifold Structure of the Weighted Projective Space.}

Throughout this paper, $a(\bf{p})$ for $a\in \bf{C}$ and ${\bf
p}=(p_{1},\cdots,p_{n})\in {\bf Z}^{n}$ will denote the diagonal
matrix:
$$Diag(a^{p_{1}},\cdots,a^{p_{n}})$$
with diagonal entries $a^{p_{i}}, i=1,\cdots,n$.  Moreover, for an
integer $q$, $\mu_{q}$ will denote the group ${\bf Z}/q\bf{Z}$.

\begin{defn}([14]) Let $Q=(q_{0},\cdots,q_{n})$ be
a $(n+1)$-tuple of positive integers. The weighted projective
space of type $Q$,
$\mbox{\bf{P}}^{n}(Q)=\mbox{\bf{P}}^{n}_{q_{0},\cdots,q_{n}}$ is
defined by

$${\bf P}^{n}_{q_{0},\cdots,q_{n}}=\left\{z\in (\mbox{\bf{C}}^{n+1})^{*}| z\sim
\lambda (\mbox{\bf{q}})\cdot z, \lambda\in \bf{C}^{*}\right\}$$
where $\lambda
(\mbox{\bf{q}})=Diag(\lambda^{q_{0}},\cdots,\lambda^{q_{n}})$.
\end{defn}

\begin{rem}
(1)The above $\bf{C}^{*}$-action is free iff $q_{i}=1$ for every
$i=0,\cdots,n$; (2) If $gcd(q_{0},\cdots,q_{n})=d\neq 1$, then
$\mbox{\bf{P}}^{n}_{q_{0},\cdots,q_{n}}$ is homeomorphic to
$\mbox{\bf{P}}^{n}_{q_{0}/d,\cdots,q_{n}/d}$(by identification of
$\lambda^{d}$ with $\lambda$ ).
\end{rem}

Weighted projective spaces are, in general, orbifolds where the
singularities have cyclic structure groups acting diagonally.
Moreover, if all the $q_{i}'s$ are mutually prime, all these
orbifold singularities are isolated. In fact, as is usually done
for complex projective spaces, we can consider the sets
$$U_{i}=\left\{[z]_{Q}\in \mbox{\bf{P}}^{n}_{q_{0},\cdots,q_{n}}: z_{i}\neq 0\right\}\subset
\mbox{\bf{P}}^{n}_{q_{0},\cdots,q_{n}}$$ and the bijective maps
$\phi_{i}$ from $U_{i}$ to ${\bf C}^{n}/\mu_{q_{i}}(Q_{i})$ given
by

$$\phi_{i}([z]_{Q})=\left(\frac{z_{0}}{(z_{i})^{q_{0}/q_{i}}},\cdots,
\frac{\hat{z_{i}}}{z_{i}},\cdots,\frac{z_{n}}{(z_{i})^{q_{n}/q_{i}}}\right)_{q_{i}}$$

where $(z_{i})^{1/q_{i}}$ is a $q_{i}$-root of $z_{i}$ and
$(\cdot)_{q_{i}}$ is a $\mu_{i}$-conjugacy class in ${\bf
C}^{n}/\mu_{q_{i}}(Q_{i})$ with $\mu_{q_{i}}$ acting on ${\bf
C}^{n}$ by
$$\xi\cdot z=\xi(Q_{i})z, \xi\in \mu_{q_{i}}$$
Here $Q_{i}=(q_{0},\cdots,\hat{q_{i}},\cdots,q_{n})$. Then on
$\phi_{i}(U_{j}\cap U_{i})\subset {\bf C}^{n}/\mu_{q_{i}}(Q_{i})$,

$$\phi_{j}\circ \phi_{i}^{-1}((z_{1},\cdots,z_{n})_{q_{i}})=
\left(\frac{z_{0}}{(z_{j})^{q_{0}/q_{j}}},\cdots,
\frac{\hat{z_{j}}}{z_{j}},\cdots,\frac{1}{(z_{j})^{q_{i}/q_{j}}},
\cdots,\frac{z_{n}}{(z_{j})^{q_{n}/q_{j}}}\right)_{q_{j}}$$

\subsection{Toric Structure of the Weighted Projective Spaces.}
Given  $Q=(q_{0},\cdots,\\ q_{n}) \in {\bf Z}^{n+1}$, let
$d_{i_{1},\cdots,i_{s}}=gcd(q_{i_{1}},\cdots,q_{i_{s}})$ and
$d_{j}=gcd(q_{0},\cdots,\hat{q}_{j},\cdots,q_{n})$ ($i_{1},\cdots,
\\ i_{s},j\in\{1,\cdots,n\}$). Define a grading of
$\mbox{\bf{C}}[X_{0},\cdots,X_{n}]$ by $degX_{i}=q_{i}$. We denote
this ring by $S(Q)$. Then ${\bf
P}^{n}_{q_{0},\cdots,q_{n}}=projS(Q)$ is the weighted projective
space of type $Q$. ${\bf P}^{n}_{q_{0},\cdots,q_{n}}$ is covered
by the affine open sets $D_{+}(X_{i}):=specS(Q)_{X_{i}},
(i=0,...,n)$. The monic monomials of $S(Q)_{X_{i}}$ are of type
$X_{i}^{-l}\Pi_{j\neq i}X_{j}^{\lambda_{j}}$, where
$lq_{i}=\Sigma_{j\neq i}\lambda_{j}q_{j}$ and $l,\lambda_{j}$ are
non-negative integers. So each such monomial is uniquely
determined by the $n$-tuple
$(\lambda_{0},\cdots,\lambda_{i-1},\lambda_{i+1},\cdots,\lambda_{n})$
of its non-negative exponents. The exponents occurring are just
the points lying in the intersection of the cone
$e:=pos\{e_{1},\cdots,e_{n}\}$ and the lattice
$N_{Q,q_{i}}\subseteq {\bf Z}^{n}$ that is defined as follows.

Consider $Q_{i}=(q_{0},\cdots,\hat{q_{i}},\cdots,q_{n})$ as an
element of $Hom_{{\bf Z}}({\bf Z}^{n},{\bf Z})$ by setting:
$$Q_{i}(a_{1},\cdots,a_{n}):=q_{0}a_{1}+\cdots,+q_{n}a_{n}$$

${\bf Z}^{n}$ being equipped with its canonical basis. Let
$\pi_{i}: {\bf Z}\longrightarrow {\bf Z}_{q_{i}}$ denote the
canonical projection. Then

$$N_{Q,q_{i}}:=Ker({\bf Z}^{n}\longrightarrow {\bf Z}\longrightarrow {\bf Z}_{q_{i}})$$
is a sublattice of ${\bf Z}^{n}$. Denote by $M_{Q,i}$ the dual
lattice. We have an isomorphism of semigroup rings
$$S(Q)_{(X_{i})}\cong \mbox{\bf{C}}[e\cap N_{Q,q_{i}}]$$
revealing $D_{+}(X_{i})$ to be the affine toric variety associated
with $\check{e}$ with respect to $M_{Q,q_{i}}$.

\begin{prop}([5]) Let
$C_{i}=(c_{1}^{i},\cdots,c_{n}^{i})$ be a basis of $N_{Q,q_{i}}$
and denote by $r_{1}^{i},\cdots,r_{n}^{i}$ the row vectors of
$C_{i}$. Let $\sigma_{i}:=pos\{r_{1}^{i},\cdots,r_{n}^{i}\}$. Then
there is an isomorphism of semigroups
$$\check{\sigma}_{i}\cap \mbox{\bf{Z}}^{n}\simeq e\cap N_{Q,q_{i}}.$$
\end{prop}

\begin{prop} ([5]) With the natation introduced
above, the matrix

$$
C_{0}= \left(
\begin{array}{cccc}
\frac{q_{0}}{d_{01}}&c_{12}^{0}&\cdots&c_{1n}^{0}\\
0&\frac{d_{01}}{d_{012}}&\cdots&c_{2n}^{0}\\
\vdots&\vdots&\ddots&\vdots\\
0&0&\cdots&\frac{d_{n}}{d_{0\cdots n}}
\end{array}
\right)= \left(
\begin{array}{c}
v_{1}\\
v_{2}\\
\vdots\\
v_{n}
\end{array}
\right) =(c_{ij}^{0})
$$
is a basis of $N_{Q,q_{0}}$, where $c_{ij}^{0}\in \bf{Z}_{\geq 0}$
are determined as follows. For fixed $j\in\{2,...,n\}$ construct
$c_{ij}^{0}$ successively for $i=j-1,...,1$ by requiring
$c_{ij}^{0}\in \bf{Z}_{\geq 0}$ to be minimal with the property
that
$$c_{ij}^{0}q_{i}+\sum_{\nu=i+1}^{j}c_{\nu j}^{0}q_{\nu}\in
 gcd(q_{0},...,q_{i-1})\mbox{\bf{Z}}=d_{0\cdots i-1}\bf{Z}.$$
\end{prop}

\begin{prop} ([5]) With the notations
introduced above, let
$v_{0}:=-\sum_{i=1}^{n}\frac{q_{i}}{q_{0}}v_{i}$,
$\rho_{i}=\mbox{\bf{R}}_{\geq 0}v_{i}$.  Then the complete fan
$\Xi$ determined by $\Xi[1]=\{\rho_{0},\cdots,\rho_{n}\}$ is the
fan of the weighted  projective space
$\mbox{\bf{P}}^{n}_{q_{0},\cdots,q_{n}}$.
\end{prop}

\begin{rem} (1) The weighted projective space
$\mbox{\bf{P}}^{n}_{q_{0},\cdots,q_{n}}$  can also be constructed
as follows. Given a fan $\Xi=\{u_{0},\cdots,u_{n}\}$ so that
$q_{0}u_{0}+q_{1}u_{1}+\cdots+q_{n}u_{n}=0$, then the toric
variety $X_{\Xi}$ is the weighted projective space ${\bf
P}^{n}_{q_{0},\cdots,q_{n}}$. From above we can see that
$\{v_{0},\cdots,v_{n}\}$ satisfies the condition
$q_{0}v_{0}+q_{1}v_{1}+\cdots+q_{n}v_{n}=0$. In fact, the
proposition 3.2.2 and 3.2.3 gave a method to compute the fan of
weighted projective space of type $Q$. (2) If
$gcd(q_{0},\cdots,q_{n})=d\neq 1$, from the construction of the
matrix $C_{0}$, we can see that
$\mbox{\bf{P}}^{n}_{q_{0},\cdots,q_{n}}$ and
$\mbox{\bf{P}}^{n}_{q_{0}/d,\cdots,q_{n}/d}$ have the same fans,
so they are homeomorphic.
\end{rem}

\subsection{Homogeneous Coordinate Representations of
Weighted Projective  Spaces Induced from Toric Varieties.}

In this section we use the theorem of David Cox (see [8]) to
represent the weighted projective space
$\mbox{\bf{P}}^{n}_{q_{0},\cdots,q_{n}}$ as the geometric quotient
$(\bf{C}^{n+1})^{*}/\bf{C}^{*}$. From the theorem of David Cox in
[7], for a fan $\Xi$, the toric variety $X_{\Xi}$ is a geometric
quotient $({\bf C})^{r}\backslash Z/G$ iff $\Xi$ is simplicial,
where $r$ is the number of the 1-primitive generators,$Z$ is a
subvariety of ${\bf C}^{r}$, $G$ is some subgroup of $({\bf
C}^{*})^{r}$. To give this representation for the weighted
projective space, we must compute the space $Z$ and the group $G$.
Let $\Xi=\{v_{0},\cdots,v_{n}\}$ be the fan of weighted projective
space $\mbox{\bf{P}}^{n}_{q_{0},\cdots,q_{n}}$. Then there are
$n+1$ 1-dimensional primitive generators, which give variables
$x_{0},\cdots,x_{n}$. Furthermore, the maximal cones of the fan
are generated  by the $n$-element subsets of
$\{v_{0},\cdots,v_{n}\}$. It follows from [9] that
$$Z=V(x_{0},\cdots,x_{n})=\{(0,\cdots,0)\}\subset {\bf C}^{n+1}.$$
Now we describe the group $G$. From [9],
$$G=\left\{(\mu_{0},\cdots,\mu_{n})\in ({\bf C}^{*})^{n+1}|
\Pi_{i=1}^{n+1}\mu_{i}^{<m,v_{i}>}=1, \mbox{for all}
 ~m\in{\bf Z}^{n}\right\}   $$
However it suffices to let $m$ be the standard basis elements
$e_{1},\cdots,e_{n}$. Thus $(\mu_{0},\cdots,\mu_{n})\in G$ iff
$$\Pi_{i=0}^{n}\mu_{i}^{<e_{1},r_{i}>}=\Pi_{i=0}^{n}\mu_{i}^{<e_{2},r_{i}>}=
\cdots=\Pi_{i=0}^{n}\mu_{i}^{<e_{n},r_{i}>}=1  \eqno{(3.1)}$$

From proposition 3.2.3, we have the vectors $v_{1},\cdots,v_{n}$,
and $v_{0}=-\Sigma_{i=1}^{n}\frac{q_{i}}{q_{0}}r_{i}$,  so we
have:

\begin{eqnarray}
\mu_{0}^{-\frac{q_{1}}{d_{01}}}\mu_{1}^{\frac{q_{0}}{d_{01}}}&=&
\mu_{0}^{(-\frac{q_{1}}{q_{0}}c_{12}^{0}-\frac{q_{2}d_{01}}{q_{0}d_{012}})}
\mu_{1}^{c_{12}^{0}}\mu_{2}^{\frac{d_{01}}{d_{012}}}=\cdots=\nonumber\\
&=&\mu_{0}^{(-\frac{q_{1}}{q_{0}}c_{1n}^{0}-\frac{q_{2}}{q_{0}}c_{2n}^{0}-\cdots-
\frac{q_{n}d_{n}}{q_{0}d_{0\cdots n}})}\mu_{1}^{c_{1n}^{0}}\cdots
 \mu_{n-1}^{c_{n-1,n}^{0}}
 \mu_{n}^{\frac{d_{n}}{d_{0\cdots n}}}=1  \nonumber
 \end{eqnarray}

So
$$
\left\{ \begin{array}{l}
\mu_{0}^{-\frac{q_{1}}{d_{01}}}=\mu_{1}^{\frac{q_{0}}{d_{01}}},\\
\mu_{0}^{(\frac{q_{1}}{q_{0}}c_{12}^{0}+\frac{q_{2}d_{01}}{q_{0}d_{012}})}
=\mu_{1}^{c_{12}^{0}}\mu_{2}^{\frac{d_{01}}{d_{012}}};\\
\cdots \cdots  \cdots \cdots \cdots \cdots \\
\mu_{0}^{(\frac{q_{1}}{q_{0}}c_{1n}^{0}+\frac{q_{2}}{q_{0}}c_{2n}^{0}+\cdots+
\frac{q_{n}d_{n}}{q_{0}d_{0\cdots n}})}=\mu_{1}^{c_{1n}^{0}}\cdots
 \mu_{n-1}^{c_{n-1,n}^{0}}
 \mu_{n}^{\frac{d_{n}}{d_{0\cdots n}}}
\end{array}
\right.
$$
We have
$$\mu_{0}=(\mu_{0}^{\frac{1}{q_{0}}})^{q_{0}},\mu_{1}=(\mu_{0}^{\frac{1}{q_{0}}})^{q_{1}},
\mu_{2}=(\mu_{0}^{\frac{1}{q_{0}}})^{q_{2}},\cdots,\mu_{n}=(\mu_{0}^{\frac{1}{q_{0}}})^{q_{n}}.
\eqno{(3.2)}$$

We obtain the following results. if $(\mu_{0},\cdots,\mu_{n})\in
G\subset({\bf C}^{*})^{n+1}$, then $(\mu_{0},\cdots,$ $\mu_{n})$
satisfies the condition (3.2). So $G\cong \bf{C}^{*}$, let
$\mu_{0}^{\frac{1}{q_{0}}}=\lambda$, then the $G$ action on
$(\bf{C}^{*})^{n+1}$ by
$\lambda(x_{0},\cdots,x_{n})=(\lambda^{q_{0}}x_{0},\cdots,\lambda^{q_{n}}x_{n})$.
Thus $\mbox{\bf{P}}^{n}_{q_{0},\cdots,q_{n}}=({\bf
C}^{n+1})^{*}/{\bf C}^{*}$.

\begin{prop} Let  $Q=(q_{0},\cdots,q_{n})$ and
$gcd(q_{0},\cdots,q_{n})=1$, then the definitions of
$\mbox{\bf{P}}^{n}_{q_{0},\cdots,q_{n}}$ in the sections 3.1 and
3.2 are equivalent.
\end{prop}


\section{The Chen-Ruan Cohomology Groups of Weighted Projective Spaces}
\subsection{The Ordinary Cohomology Groups of Weighted
Projective Spaces.}  Let $Q=(q_{0},\cdots,q_{n})$,
$gcd(q_{0},\cdots,q_{n})=1$ and
$\mbox{\bf{P}}^{n}_{q_{0},\cdots,q_{n}}$  be the weighted
projective space of the type $Q$. The ordinary cohomology group of
$\mbox{\bf{P}}^{n}_{q_{0},\cdots,q_{n}}$ has already been studied
by several authors,  [1], [17]. Here we only give the results.

\begin{thm} Let $Q=(q_{0},\cdots,q_{n})$ and
$\mbox{\bf{P}}^{n}_{q_{0},\cdots,q_{n}}$  be the weighted
projective space of the type $Q$, then the cohomology group of
$\mbox{\bf{P}}^{n}_{q_{0},\cdots,q_{n}}$ with rational coefficient
is

$$
H^{i}(\mbox{\bf{P}}^{n}_{q_{0},\cdots,q_{n}},\mbox{\bf{Q}})=\left\{
\begin{array}{ll}
\bf{Q},& \mbox{if}~~  i=2r,0\leq r\leq n;\\
0,& \mbox{if}~~~ i ~~~ \mbox{is odd or}~~ i>2n.
\end{array}
\right.
$$
\end{thm}

\subsection{Orbiford Structure from Toric Varieties.}
Let $Q=(q_{0},\cdots,q_{n})$ and $\mbox{\bf{P}}^{n}_{q_{0}
\cdots,q_{n}}$ be the weighted projective space of the type $Q$.
And let $\Xi=\{v_{0},\cdots,v_{n}\}$ be the fan of
$\mbox{\bf{P}}^{n}_{q_{0},\cdots,q_{n}}$ in the lattice ${\bf
Z}^{n}$ defined in section 3.2. We have the following proposition.

\begin{prop} Let $Q=(q_{0},\cdots,q_{n})$ and
$\mbox{\bf{P}}^{n}_{q_{0},\cdots,q_{n}}$ be the weighted
projective space of the type $Q$. Then taking
$Q=(q_{0},\cdots,q_{n})$ as toric variety and its fan $\Xi$ is
generated by $\{v_{0},\cdots,v_{n}\}$ above, it has the orbifold
structure as follows.
$$\{(U_{\sigma_{k}^{'}},G_{\sigma_{k}},\pi_{\sigma_{k}})| \sigma_{k}=(v_{0},
\cdots,\hat{v}_{k},\cdots,v_{n}), k=0,...,n \}$$ In particular,
$G_{\sigma_{k}}={\bf Z}_{q_{k}}$ is the cyclic group of order
$q_{k}$.
\end{prop}

\begin{proof} Because the fan $\Xi$ of toric variety
$\mbox{\bf{P}}^{n}_{q_{0},\cdots,q_{n}}$ is generated by
$\{v_{0},\cdots,v_{n}\}$, so we all have $n+1$ $n$-dimensional
cones $\sigma_{0},\cdots,\sigma_{n}$. From section 2.4, we can
conclude that
$$\{(U_{\sigma_{k}^{'}},G_{\sigma_{k}},\pi_{\sigma_{k}})|
\sigma_{k}=(v_{0}, \cdots,\hat{v}_{k},\cdots,v_{n}), k=0,...,n
\}$$forms the orbifold structure of
$\mbox{\bf{P}}^{n}_{q_{0},\cdots,q_{n}}$. Now we prove that for
any $\sigma_{k}$, $G_{\sigma_{k}}={\bf Z}_{q_{k}}$. From section
2.4 we have $G_{\sigma_{k}}=N/N_{\sigma_{k}}$, let $\pi_{k}:{\bf
Z}\longrightarrow {\bf Z}_{q_{k}}$ be the standard projection.
Define the map $Q_{k}:{\bf Z}^{n}\longrightarrow {\bf Z}$ such
that $Q_{k}(a_{1},\cdots,a_{n})=q_{0}a_{1}+\cdots+q_{n}a_{n},
Q_{i}=(q_{0},\cdots,\hat{q_{k}},\cdots,q_{n})$. Then we let
$$N_{Q,q_{k}}:=Ker({\bf Z}^{n}\longrightarrow {\bf Z}\longrightarrow {\bf Z}_{q_{k}})$$
From H.Conrads [5], because $\Xi=\{v_{0},\cdots,v_{n}\}$ is the
fan of $\mbox{\bf{P}}^{n}_{q_{0},\cdots,q_{n}}$, so $N_{Q,q_{k}}$
is generated by
$\sigma_{k}=(v_{0},\cdots,\hat{v_{k}},\cdots,v_{n})$, i.e.
$N_{Q,k}=N_{\sigma_{k}}$. So we have:
$$G_{\sigma_{k}}=N/N_{\sigma_{k}}={\bf Z}^{n}/N_{Q,q_{k}}\cong {\bf Z}_{q_{k}}.$$
\end{proof}

\begin{rem} From the proposition above,
$U_{\sigma_{k}^{'}}={\bf C}^{n}$, so
$U_{\sigma_{k}}=U_{\sigma_{k}^{'}}/G_{\sigma_{k}}={\bf C}^{n}/{\bf
Z}_{q_{k}}$, and the action of $\mbox{\bf{Z}}_{q_{k}}$ is the
diagonal action, its  matrix representation can be computed from
the method of section 2.4 which we will use in the following
section.
\end{rem}

\subsection{Twisted Sectors of Weighted Projective Spaces
and Degree Shifting Numbers.} Given a weighted projective space
${\bf P}(Q)$ of type $Q=(q_{0},\cdots,q_{n})$, we now discuss the
twisted sectors of ${\bf P}(Q)$. First, taken as toric variety,
${\bf P}(Q)$ has the orbit decomposition ${\bf
P}(Q)=\bigsqcup_{\tau \in \Xi}O_{\tau}$, we will determined
$O_{\tau}$. From section 2.4, if $\sigma\in \Xi$ is a
$n$-dimensional cone, suppose $\sigma=(v_{0},\cdots,v_{n-1})$,
then $G_{\sigma}={\bf Z}_{q_{n}}$. We also know that
$G_{\sigma}=\{k_{a}=\Sigma_{i=0}^{n-1}a_{i}v_{i}: k_{a}\in N,
a_{i}\in [0,1)\}$. If $\tau=(v_{0},\cdots,v_{i-1})$ is  a face of
$\sigma$, then $G_{\tau}=\{g_{a}\in G_{\sigma}: a_{i}=0 ~\mbox{if}
~j+1\leq i\leq n \}$, i.e.
$G_{\tau}=\{k_{a}=\Sigma_{j=1}^{i-1}a_{j}v_{j}: k_{a}\in
N,a_{j}\in[0,1)\}$. $G_{\tau}$ can be taken as the local group of
the points in $O_{\tau}$. For the weighted projective space, we
have the following proposition.

\begin{prop} Given a weighted projective space
${\bf P}(Q)$ of type $Q=(q_{0},\cdots,q_{n})$. Let
$\Xi=\{v_{0},\cdots,v_{n}\}$ be the fan of ${\bf P}(Q)$. If
$\tau=(v_{0},\cdots,v_{i-1})$is a cone of $\Xi$, then we have
$G_{\tau}\cong {\bf Z}_{d}$, where $d=gcd(q_{i},\cdots,q_{n})$. In
particular, $(q_{i},\cdots,q_{n})$ is the maximal subset of
$(q_{0},\cdots,q_{n})$ whose $gcd$ is $d$ iff the dimension of
fixed point set of ${\bf Z}_{d}$ is $n-i$.
\end{prop}

\begin{proof} Let
$\sigma_{k}=(v_{0},\cdots,v_{i-1},v_{i},\cdots,\hat{v_{k}},\cdots,v_{n}),
i\leq k\leq n$, then from proposition 4.2.1, $G_{\sigma_{k}}\cong
{\bf Z}_{q_{k}}$. Since we have
$G_{\sigma_{k}}=\{k_{a}=\Sigma_{j\neq k}a_{j}v_{j}: k_{a}\in N,
a_{j}\in [0,1),i\leq k\leq n\}$,
$G_{\tau}=\{k_{a}=\Sigma_{j=1}^{i-1}a_{j}v_{j}: k_{a}\in N,
a_{j}\in [0,1)\}$, so $G_{\tau}$ is a subgroup of
$G_{\sigma_{k}}(i\leq k\leq n)$. Next we show that any element of
a common subgroup of $G_{\sigma_{k}}(i\leq k\leq n)$ has the form
of the elements in $G_{\tau}$. For any element $g_{a}\in {\bf
Z}_{q_{k}}(i\leq k\leq n)$. Without loss of generality, assume
$g_{a}\in {\bf Z}_{q_{i}},{\bf Z}_{q_{i+1}}$, and
$$(1)g_{a}=a_{0}v_{0}+\cdots+a_{i-1}v_{i-1}+a_{i+1}v_{i+1}+\cdots+a_{n}v_{n}$$
$$(2)g_{a}=a_{0}^{'}v_{0}+\cdots+a_{i}^{'}v_{i}+a_{i+2}^{'}v_{i+2}+\cdots+a_{n}^{'}v_{n}$$
Since the action of the same $g_{a}$ at one neighborhood of a
point is the same, so for $j\neq i,i+1, a_{j}=a_{j}^{'}$,
$(1)-(2)$, we obtain $a_{i+1}v_{i+1}-a_{i}^{'}v_{i}=0$, and we
have a contradiction, because the generators $v_{i+1}$ and $v_{i}$
of the fan are independent over $N_{{\bf R}}$. We have
$g_{a}=a_{0}v_{0}+\cdots+a_{i-1}v_{i-1}$ and $g_{a}\in G_{\tau}$.
So we prove $G_{\tau}={\bf Z}_{d}, d=gcd(q_{i},\cdots,q_{n})$.

Now if $(q_{i},\cdots,q_{n})$ is the maximal subset of
$(q_{0},\cdots,q_{n})$ that satisfies the condition
$gcd(q_{i},\cdots,q_{n})=d$, while the dimension of fixed point
set of ${\bf Z}_{d}$ is not $n-i$. Because the dimension of
$O_{\tau}$ is $n-i$, then we have a generator $g_{a}\in {\bf
Z}_{d}$ and $g_{a}$ can be represented by
$g_{a}=a_{0}v_{0}+\cdots+a_{i-1}v_{i-1}$. And we must have some
$a_{s}=0$ for $s\leq i-1$. If we let
$\rho=(v_{0},\cdots,\hat{v_{s}},\cdots,v_{i-1})$, then from the
first part of the proposition we have $G_{\rho}={\bf
Z}_{d}=gcd(q_{s},q_{i},\cdots,q_{n})$ which will contradict the
maximal principle of the condition.

Conversely, suppose the dimension of the fixed point set of ${\bf
Z}_{d}$ is $n-i$, i.e. the dimension of the orbit $O_{\tau}$. If
we have a subset of $(q_{0},\cdots,q_{n})$ whose great common
divisor is $d$, and the number of the subset is more than the
$(q_{i},\cdots,q_{n})$. Without loss of generality, we assume
$gcd(q_{s},q_{i},\cdots,q_{n})=d$ for $0\leq s\leq i-1$. Then let
$\delta=(v_{0},\cdots,v_{s-1},v_{s+1},\cdots,v_{i-1})$.  From the
first part of the theorem we have
$$G_{\delta}={\bf Z}_{d}=\{k_{a}=a_{0}v_{0}+\cdots+a_{s-1}v_{s-1}
+a_{s+1}v_{s+1}+\cdots+a_{i-1}v_{i-1}: a_{i}\in [0,1)\}$$ We see
that the dimension of the fixed point set of ${\bf Z}_{d}$
overdues $n-i+1$, a contradiction.
\end{proof}

Now we discuss the twisted  sectors of the weighted projective
space ${\bf P}(Q)$. From the theorem of Poddar in [21] for the
twisted sectors of general toric varieties, we have:

\begin{thm}([21]) A twisted sector of a weighted
projective space is isomorphic to a subvariety
$\overline{O}_{\tau}$ of $X_{\Xi}={\bf P}(Q)$ for some $\tau \in
\Xi$. There is a one to one  correspondence between the set of
twisted sectors of the type $\overline{O}_{\tau}$ and the set of
integral vectors in the interior of $\tau$ which are linear
combinations of the 1-dimensional generators of $\tau$ with
coefficients in $(0,1)$.
\end{thm}

For the orbit $O_{\tau}$ of the weighted projective space ${\bf
P}(Q)$ for some $\tau \in \Xi$, we have the following proposition:

\begin{prop} Given a weighted projective space
${\bf P}(Q)$ of type $Q=(q_{0},\cdots,q_{n})$. Let
$\Xi=\{v_{0},\cdots,v_{n}\}$ be the fan of ${\bf P}(Q)$. If
$\tau=(v_{0},\cdots,v_{i-1})$ is a cone in $\Xi$, then
$\overline{O}_{\tau}$ is the weighted projective space ${\bf
P}(Q_{\tau})$, where $Q_{\tau}=(q_{i},\cdots,q_{n})$.
\end{prop}

\begin{proof} From Fulton [13], $\overline{O}_{\tau}$ is a
toric variety and the fan $Star(\tau)$ can be described as
follows. Let $N_{\tau}$ be the sublattice of $N$ generated by
$\tau\in N$, and $N(\tau)=N/N_{\tau}, M(\tau)=\tau^{\perp}\cap M$
be the quotient lattice and the dual. The star of a cone $\tau$
can be defined abstractly as the set of cones $\sigma$ in $\Xi$
that contain $\tau$ as a face. Such cones $\sigma$ are determined
by their images in $N(\tau)$, i.e. by
$$\overline{\sigma}=(\sigma + (N_{\tau})_{{\bf R}})/(N_{\tau})_{{\bf R}}\subset
N_{{\bf R}}/(N_{\tau})_{{\bf R}}=N(\tau)_{{\bf R}}$$ These cones
$\{\overline{\sigma}: \tau< \sigma\}$ form a fan in $N(\tau)$, and
we denote this fan by $Star(\tau)$, the corresponding toric
variety is $n-k$-dimensional toric variety. For the toric variety
${\bf P}(Q)$, let $\overline{v}_{i},\cdots,\overline{v}_{n}$ be
the images of $v_{0},\cdots,v_{n}$ in the quotient lattice
$N(\tau)$, since in $N$, $q_{0}v_{0}+\cdots+q_{n}v_{n}=0$, we have
$q_{i}\overline{v}_{i}+\cdots+q_{n}\overline{v}_{n}=0$ in the
quotient lattice $N(\tau)$. So from the definition of weighted
projective space, Remark 3.2.4, we  conclude  that the toric
variety corresponding to the fan $Star(\tau)$ in the quotient
lattice $N(\tau)$ is the weighted projective space ${\bf
P}(Q_{\tau})$, where $Q_{\tau}=(q_{i},\cdots,q_{n})$.
\end{proof}

\begin{rem} The $gcd$ of $(q_{i},\cdots,q_{n})$
need not necessary be 1, if $d=gcd(q_{i},\cdots,q_{n})$, in
general, ${\bf P}(Q_{\tau})$ is a nonreduced orbifold, and it has
a corresponding reduced orbifold  ${\bf P}(Q_{\tau}^{'})$, where
$Q_{\tau}^{'}=(q_{i}/d,\cdots,q_{n}/d)$.
\end{rem}

\begin{thm} Let $Q=(q_{0},\cdots,q_{n})$ and  ${\bf
P}(Q)= {\bf P}^{n}_{q_{0},\cdots,q_{n}}$ be the weighted
projective space of the type $Q$. Let $\Xi=\{v_{0},\cdots,v_{n}\}$
be the fan of $\mbox{\bf{P}}(Q)$. If $\tau=(v_{0},\cdots,v_{i-1})$
is a cone in $\Xi$, then $\overline{O}_{\tau}={\bf P}(Q_{\tau})$,
where $Q_{\tau}=(q_{i},\cdots,q_{n})$. Let
$gcd(q_{i},\cdots,q_{n})=d$. Then we have:
$\overline{O}_{\tau}={\bf P}(Q_{\tau})$ is a twisted sector iff
$(q_{i},\cdots,q_{n})$ is the maximal subset of
$(q_{0},\cdots,q_{n})$ that satisfies the condition
$gcd(q_{i},\cdots,q_{n})=d$.
\end{thm}

\begin{proof} If $\overline{O}_{\tau}$ is a twisted sector,
then suppose $\overline{O}_{\tau}=X_{(g_{a})}$, where $g_{a}\in
G_{\tau}={\bf Z}_{d}$ is a generator. If we have
$q_{k}\bar{\in}(q_{i},\cdots,q_{n})$,
$gcd(q_{k},q_{i},\cdots,q_{n})=d$, let
$\delta=(v_{0},\cdots,\hat{v_{k}},\cdots,v_{i-1})$, then
$\overline{O}_{\delta}={\bf P}(Q_{\delta})$, where
$Q_{\delta}=(q_{k},q_{i},\cdots,q_{n})$. Since from section 2.4,
$G_{\delta}={\bf Z}_{d}=\{k_{a}=\Sigma_{j=0}^{i-1}a_{j}v_{j}:
j\neq k, a_{j}\in [0,1)\}$,  the coefficients of the
representation of $g_{a}$ are non zeroes, so from theorem 4.3.2,
$X_{(g_{a})}=\overline{O}_{\delta}$ which contradicts the
condition.

Suppose $(q_{i},\cdots,q_{n})$ is the maximal subset of
$(q_{0},\cdots,q_{n})$ that meets the condition
$gcd(q_{i},\cdots,q_{n})=d$,  from proposition 4.4, let
$\tau=(v_{0},\cdots,v_{i-1}), G_{\tau}\cong {\bf Z}_{d}$ and the
dimension of the fixed point set of ${\bf Z}_{d}$ is $n-i$, the
dimension of the orbit $O_{\tau}$. Since
$$G_{\tau}={\bf
Z}_{d}=\{k_{a}=\Sigma_{j=0}^{i-1}a_{j}v_{j}: a_{j}\in [0,1)\}$$ we
must have one generator $g_{a}\in G_{\tau}$ so that $g_{a}$ can be
represented by
$$g_{a}=a_{0}v_{0}+\cdots+a_{i-1}v_{i-1}, \mbox{all}~~ a_{j}\neq 0, j=0,\cdots,i-1.$$
Thus the dimension of the fixed point set of $g_{a}$ is $n-i$.
From theorem 4.3.2, we  see that $X_{(g_{a})}=\bar{O}_{\tau}$.
\end{proof}

\begin{rem} From above analysis, if
$\tau=(v_{0},\cdots,v_{i-1})$ is a cone in $\Xi$, we describe the
orbifold structure of twisted sector $X_{(g_{a})}=\bar{O}_{\tau}$
as follows. In the points of $O_{\tau}$, the local group at the
point $x$ is ${\bf Z}_{d}$, $d=gcd(q_{i},\cdots,q_{n})$, and
$C(g_{a})={\bf Z}_{d}$ acts trivially. If $\delta>\tau$ is a cone,
then $O_{\delta}\subset \overline{O}_{\tau}$,the local group at
the $y\in O_{\delta}$ is the cyclic group ${\bf Z}_{t}$, where
$t=gcd(q_{t_{1}},\cdots,q_{t_{s}}),q_{i},\cdots,q_{n}$, and
$\delta=(v_{0},\cdots,\hat{v}_{t_{1}},\cdots,\hat{v}_{t_{s}},\cdots,v_{i-1})$.
And the action of ${\bf Z}_{t}$ on $\bar{O}_{\tau}$ can be
described as before.
\end{rem}

The degree shifting numbers can be computed easily. For instance,
let $X_{(g_{a})}=\overline{O}_{\tau}$ be a twisted sector, and
$\tau=(v_{0},\cdots,v_{i-1})$, we can write $g_{a}$ as
$$g_{a}=\Sigma_{j=0}^{i-1}a_{j}v_{j},  a_{j}\in (0,1)$$
The degree shifting number $\iota_{(g_{a})}$ of $X_{(g_{a})}$ is:
$$\iota_{(g_{a})}=\Sigma_{j=0}^{i-1}a_{j}$$
So in the orbifold structure of weighted projective space ${\bf
P}(Q)$, we can compute the degree shifting number corresponding to
any twisted sector.
\subsection{The Chen-Ruan Cohomology Groups of Weighted
Projective Spaces.}

Up to now, given a weighted projective space ${\bf P}(Q)$ of type
$Q=(q_{0},\cdots,q_{n})$, from Theorem 4.3.5 we can calculate the
twisted sectors of ${\bf P}(Q)$ and we also can compute the degree
shifting  numbers of the corresponding twisted sectors. So  we
write the Chen-Ruan  cohomology group of ${\bf P}(Q)$  in the
following manner.

\begin{thm} Let $Q=(q_{0},\cdots,q_{n})$ and ${\bf
P}(Q)= \mbox{\bf{P}}^{n}_{q_{0},\cdots,q_{n}}$ be the weighted
projective space of the type $Q$.  Let
$\Xi=\{v_{0},\cdots,v_{n}\}$ be the fan of $\mbox{\bf{P}}(Q)$.
Then the orbifold cohomology group of ${\bf P}(Q)$ is

$$H_{orb}^{p}({\bf P}(Q);{\bf Q})\cong \bigoplus_{\sigma\in \Xi,l\in {\bf Q}}
H^{p-2l}(\overline{O}_{\sigma})\otimes \oplus_{t\in \sigma_{l}}
{\bf Q}t$$ Where $\sigma_{l}=\{\sum_{ v_{i}\subset \sigma}
a_{i}v_{i}\in N: a_{i}\in(0,1),\sum_{v_{i}\subset
\sigma}a_{i}=l\}$.(when $\sigma=0$, set $l=0, T(\sigma)_{l}={\bf
C}$). $\overline{O}_{\sigma}$ is the closure of the orbit
corresponding to $\sigma \in\Xi$. Here $p$ is rational numbers in
$[0,n]$, and $H^{p-2l}(\overline{O}_{\sigma})=0$ if $p-2l$ is not
integral.
\end{thm}

Note that the elements of $\oplus_{0\neq \sigma\in
\Delta,l}\sigma_{l}$ correspond to the twisted sectors of ${\bf
P}(Q)$.


\subsection{Example.}  For $Q=(2,3,4)$, ${\bf P}(Q)={\bf
P}_{2,3,4}^{2}$, we have $q_{0}=2,q_{1}=3,q_{2}=4$. From
proposition 3.2.2., we have
$$
C_{0}= \left(
\begin{array}{cc}
2&0\\
0&1
\end{array}
\right)
$$
So let $v_{1}=(2,0), v_{2}=(0,1)$, we have
$v_{0}=-\frac{q_{1}}{q_{0}}v_{1}-\frac{q_{2}}{q_{0}}v_{2}=(-3,-2)$.
The fan $\Xi$ of ${\bf P}_{2,3,4}^{2}$ is generated by
$\{v_{0},v_{1},v_{2}\}$. For
$\sigma_{2}=(v_{0},v_{1})=((-3,-2),(2,0))$, we have
$G_{\sigma_{2}}=N/N_{\sigma_{2}}={\bf Z}_{4}$. We  write the
matrix representation of the action of ${\bf Z}_{4}$ on
$U_{\sigma_{2}^{'}}={\bf C}^{2}$ as follows:
$$
\left(
\begin{array}{cc}
1&0\\
0&1
\end{array}
\right); \left(
\begin{array}{cc}
e^{2\pi i\cdot \frac{1}{2}}&0\\
0&e^{2\pi i\cdot \frac{1}{4}}
\end{array}
\right) ; \left(
\begin{array}{cc}
1&0\\
0&e^{2\pi i\cdot \frac{1}{2}}
\end{array}
\right) ; \left(
\begin{array}{cc}
e^{2\pi i\cdot \frac{1}{2}}&0\\
0&e^{2\pi i\cdot \frac{3}{4}}
\end{array}
\right)
$$
For $\sigma_{1}=(v_{0},v_{2})=((-3,-2),(0,1))$, we have
$G_{\sigma_{1}}={\bf Z}_{3}$. We  write the matrix representation
of the action of ${\bf Z}_{3}$ on $U_{\sigma_{1}^{'}}={\bf C}^{2}$
as follows:
$$
\left(
\begin{array}{cc}
1&0\\
0&1
\end{array}
\right); \left(
\begin{array}{cc}
e^{2\pi i\cdot \frac{1}{3}}&0\\
0&e^{2\pi i\cdot \frac{2}{3}}
\end{array}
\right) ;  \left(
\begin{array}{cc}
e^{2\pi i\cdot \frac{2}{3}}&0\\
0&e^{2\pi i\cdot \frac{1}{3}}
\end{array}
\right)
$$
And for $\sigma_{0}=(v_{1},v_{2})=((2,0),(0,1))$, we have
$G_{\sigma_{0}}={\bf Z}_{2}$ and representation:
$$
\left(
\begin{array}{cc}
1&0\\
0&1
\end{array}
\right);  \left(
\begin{array}{cc}
e^{2\pi i\cdot \frac{1}{2}}&0\\
0&1
\end{array}
\right)
$$
If we let
$$ g_{1}=\left(
\begin{array}{cc}
e^{2\pi i\cdot \frac{1}{2}}&0\\
0&e^{2\pi i\cdot \frac{1}{4}}
\end{array}
\right); g_{2}=\left(
\begin{array}{cc}
e^{2\pi i\cdot \frac{1}{3}}&0\\
0&e^{2\pi i\cdot \frac{2}{3}}
\end{array}
\right); g_{3}=\left(
\begin{array}{cc}
e^{2\pi i\cdot \frac{1}{2}}&0\\
0&1
\end{array}
\right)
$$
then we have the twisted  sectors: let $p_{1}=[0,1,0],
p_{2}=[0,0,1]$, $X_{(g_{1})}=X_{(g_{1}^{3})}=p_{2},
\iota_{(g_{1})}=\frac{3}{4}, \iota_{(g_{1}^{3})}=\frac{5}{4}$;
$X_{(g_{1}^{2})}$ and $X_{(g_{3})}$ are the same twisted sector
${\bf P}_{2,0,4}$, and
$\iota_{(g_{1}^{2})}=\iota_{(g_{3})}=\frac{1}{2}$;
$X_{(g_{2})}=X_{(g_{2}^{2})}=p_{1},
\iota_{(g_{2})}=\iota_{(g_{2}^{2})}=1$. So
\begin{eqnarray}
H^{p}_{orb}({\bf P}_{2,3,4}^{2}; {\bf Q})&=&H^{p}({\bf
P}_{2,3,4}^{2}; {\bf Q}) \oplus H^{p-2\iota_{(g_{1})}}(\{p_{2}\};
{\bf Q})\oplus H^{p-2\iota_{(g_{1}^{3})}}(\{p_{2}\}; {\bf Q}) \nonumber \\
& \oplus & H^{p-2\iota_{(g_{1}^{2})}}({\bf P}_{2,0,4}; {\bf
Q})\oplus 2H^{p-2\iota_{(g_{2})}}(\{p_{1}\}; {\bf Q}) \nonumber
\end{eqnarray}
We compute the orbifold cohomology group of ${\bf P}^{2}_{2,3,4}$
as
$$
\begin{array}{l}
H^{0}_{orb}({\bf P}^{2}_{2,3,4}; {\bf Q})={\bf Q};\\
H^{1}_{orb}({\bf P}^{2}_{2,3,4}; {\bf Q})={\bf Q};\\
H^{\frac{3}{2}}_{orb}({\bf P}^{2}_{2,3,4}; {\bf Q})={\bf Q};\\
H^{2}_{orb}({\bf P}^{2}_{2,3,4}; {\bf Q})={\bf Q}\oplus {\bf Q}\oplus {\bf Q};\\
H^{\frac{5}{2}}_{orb}({\bf P}^{2}_{2,3,4}; {\bf Q})={\bf Q};\\
H^{3}_{orb}({\bf P}^{2}_{2,3,4}; {\bf Q})={\bf Q};\\
H^{4}_{orb}({\bf P}^{2}_{2,3,4}; {\bf Q})={\bf Q}.
\end{array}
$$
All the other dimensions of the Chen-Ruan cohomology groups are
zero.

\section{The Chen-Ruan  Cohomology Rings of the Weighted Projective Spaces}
\subsection{The Ordinary Cohomology Ring of Weighted Projective
Spaces.} In this section we recall the ordinary cohomology ring of
the weighted projective space. The readers may refer to [1]. Let
$Q=(q_{0},\cdots,q_{n})$ and ${\bf P}(Q)=
\mbox{\bf{P}}^{n}_{q_{0},\cdots,q_{n}}$ be the weighted projective
space of the type $Q$. Let ${\bf P}^{n}$ be the $n$-dimensional
complex projective space. As in [1], let $\varphi: {\bf
P}^{n}\longrightarrow {\bf P}(Q)$ be the map taking
$[x_{0},\cdots,x_{n}]$ to $[x_{0}^{q_{0}},\cdots,x_{n}^{q_{n}}]$.
Take $k\in\{0,\cdots,n\}$, and consider $I=\{i_{0},\cdots,i_{k}\}$
with $0\leq i_{0}<\cdots<i_{k}\leq n$. Put
$l_{I}=l_{I}(q_{i_{0}},\cdots,q_{i_{k}})=q_{i_{0}}\cdots
 q_{i_{k}}/gcd(q_{i_{0}},\cdots,q_{i_{k}})$, and let
$$l_{k}=l_{k}(q_{0},\cdots,q_{n})=lcm\{l_{I}|I\subset
\{0,\cdots,n\},|I|=k+1\}$$

\begin{thm}([1]) For each $k, 0\leq k\leq n$,
there exists a unique $\xi_{k}\in H^{2k}({\bf P}(Q);{\bf Q})$ such
that $\varphi^{*}(\xi_{k})=l_{k}\xi^{k}$, and
$\{1,\xi,\cdots,\xi^{n}\}$ is a {\bf Q}-basis of the free abelian
group $H^{2k}({\bf P}(Q);{\bf Q})$. In other words there are
commutative diagrams
\begin{center}
\setlength{\unitlength}{0.8cm}
\begin{picture}(5,3)\put(0,2){\shortstack{$H^{2k}({\bf P}(Q);{\bf Q})$}}
\put(3.7,2.3){\shortstack{$\varphi^{*}$}}
 \put(2.9,2.1){\vector(1,0){2.0}}
 \put(5,2){\shortstack{$H^{2k}({\bf P}(Q);{\bf Q})$}}
 \put(1.2,1.3){\shortstack{$\|$}}
 \put(6,1.3){\shortstack{$\|$}}
 \put(1.1 ,0.7){\shortstack{${\bf Q}$}}
\put(5.9,0.7){\shortstack{${\bf Q}$}}
 \put(2.9,0.8){\vector(1,0){2.0}}
 \put(3.5,0.4){\shortstack{$\cdot l_{k}$}}
\end{picture}
\end{center}
\end{thm}
So  we can make precise the  multiplicative structure of the
cohomology $H^{2k}({\bf P}(Q);{\bf Q})$. Since $\varphi^{*}:
H^{*}({\bf P}(Q);{\bf Q})\longrightarrow H^{*}({\bf P}^{n};{\bf
Q})$ is a ring homomorphism, so
$$
\xi_{i}\xi_{j}=\left\{
\begin{array}{ll}
e_{ij}\xi_{i+j},&  i+j \leq n;\\
0,& \mbox{if}~~\mbox{not}.
\end{array}
\right.
$$
Where $e_{ij}=l_{i}l_{j}/l_{i+j}, 1\leq i,j\leq n$.

In the polynomial  ring ${\bf Q}[T_{1},\cdots,T_{n}]$, let $A$ be
the idea generated by the elements
$$T_{i}T_{j}(i+j>n)~~ \mbox{and}~~ T_{i}T_{j}-e_{ij}T_{i+j}(i+j\leq n)$$
We obtain a ring isomorphism
$$H^{*}({\bf P}(Q);{\bf Q})\cong {\bf Q}[T_{1},\cdots,T_{n}]/A $$
Where $\xi_{i}$ corresponds to the class of $T_{i}$.

\subsection{Three Multi-sectors.}  Let ${\bf P}(Q)$ be
the weighted projective space of type $Q=(q_{0},\cdots,q_{n})$,
and $\Xi=\{v_{0},\cdots,v_{n}\}$ be the fan when ${\bf P}(Q)$ is
taken as toric variety. Denote $\Xi[n]$ the set of $n$-dimensional
cones. For a cone $\tau\in \Xi$, denote the set of its primitive
1-dimensional generators by $\tau[1]$. From the section 2.4,
$G_{\tau}=\{\Sigma_{v_{i}\subset \tau[1]}a_{i}v_{i}:
a_{i}\in[0,1)\}$, let $R(\tau):=\{g_{a}=\Sigma a_{i}u_{i}|u_{i}\in
\tau[1],0\leq a_{i}<1\}\cap N$. Now we describe the 3-multisector
$X_{({\bf g})}$ for ${\bf g}=(g_{1},g_{2},g_{3})$.

Take any $x\in X_{\Xi}={\bf P}(Q)$ with  nontrivial local group.
Then $x$ belongs to a unique $O_{\tau}$ such that $\tau$ is not
the trivial cone. Pick any elements $g_{a},g_{b}$ from
$G_{x}=G_{\tau}$. We shall find $X_{({\bf g})}$ where ${\bf
g}=(g_{a},g_{b},(g_{a}g_{b})^{-1})$. Let $\tau_{a}$ and $\tau_{b}$
be the faces of $\tau$, whose interiors contain $g_{a}$ and
$g_{b}$ respectively. Let $\sigma$ be any $n$-dimensional cone
containing $\tau$. Let $z$ be any point in $U_{\sigma^{'}}$.
Suppose $z$ is fixed by both $g_{a}$ and $g_{b}$. Then
$\chi^{u^{i}}=0$ whenever $u_{i}\in\tau_{a}\cup\tau_{b}$. Hence
$\pi_{\sigma}(z)\in
\overline{O}_{\tau_{a}}\cap\overline{O}_{\tau_{b}}\cap
U_{\sigma}$. A local uniformizing system for $X_{({\bf g})}$ is
given by $(V_{x}^{{\bf g}},G_{x},\pi)$, where
$$V_{x}^{{\bf g}}=({\bf C}^{j}\times W)\cap \{\chi^{u_{i}}=0: \forall
 u_{i}\in\tau_{a}\cup\tau_{b}\}$$
This leads us to observe that $\{(x,{\bf g})\in X_{({\bf g})}|x\in
U_{\sigma}\}$ is complex analytical isomorphic to
$\overline{O}_{\tau_{a}}\cap \overline{O}_{\tau_{b}}\cap
U_{\sigma}$. Since this is true in respective of the choice of
$\sigma$, $X_{({\bf g})}\cong \overline{O}_{\tau_{a}}\cap
\overline{O}_{\tau_{b}}$. So we have:

\begin{thm}([22]) If $\tau_{1}[1]\cup \tau_{2}[1]$
generates an element of $\Xi$, then for every pair $g_{a_{1}}\in
R(\tau_{1})\cap Int(\tau_{1}), g_{a_{2}}\in R(\tau_{2})\cap
Int(\tau_{2})$, we have a unique 3-multisector $X_{({\bf g})}$,
which is analytically isomorphic to $\overline{O}_{\tau_{1}}\cap
\overline{O}_{\tau_{2}}$. As we vary over $\tau_{1},\tau_{2}$, we
obtain all the 3-multisectors.
\end{thm}

Since the fan $\Xi=\{v_{0},\cdots,v_{n}\}$ of ${\bf P}(Q)$ for
$Q=(q_{0},\cdots,q_{n})$ has $n+1$ primitive 1-dimensional
generators. If $\tau_{1}$ and $\tau_{2}$ are two cones of $\Xi$,
then we have that $\tau_{1}[1]\cup \tau_{2}[1]$ form an element
$\tau=\tau_{1}[1]\cup \tau_{2}[1]$ of $\Xi$, so
$\overline{O}_{\tau_{1}}\cap \overline{O}_{\tau_{2}}$ is an
3-multisector of ${\bf P}(Q)$. Moreover, we can prove the
3-multisectors of ${\bf P}(Q)$ are actually twisted sectors.

\begin{thm} Let $X={\bf P}(Q)$ be the weighted
projective space of type $Q=(q_{0},\cdots,q_{n})$. $X_{(g_{1})}$
and $X_{(g_{2})}$ are two twisted sectors of $X$, suppose they
correspond to the cones $\tau_{1}$ and $\tau_{2}$ respectively,
i.e.
$X_{(g_{1})}=\overline{O}_{\tau_{1}}$,$X_{(g_{2})}=\overline{O}_{\tau_{2}}$.
Then $X_{({\bf g})}=X_{(g_{1},g_{2},(g_{1}g_{2})^{-1})}$ is still
a twisted sector.
\end{thm}

\begin{proof} First  if $\tau_{1}\subset \tau_{2}$, then
from theorem 5.2.1, $\tau_{1}[1]\cup \tau_{2}[1]$ generates
$\tau_{2}$, so $X_{({\bf
g})}=\overline{O}_{\tau_{2}}=X_{(g_{2})}$.

If $\tau_{1}[1]\cap \tau_{2}[1]=\emptyset$, let
$\tau=(\tau_{1}[1]\cup \tau_{2}[1])$, then $X_{({\bf
g})}=X_{(g_{1},g_{2},(g_{1}g_{2})^{-1})}=\overline{O}_{\tau}$.
Since $G_{\tau}=\{\Sigma_{v_{i}\subset \tau}a_{i}v_{i}|a_{i}\in
[0,1)\}$, we always can find an element $g\in G_{\tau}$ such that
$g=\Sigma_{v_{i}\subset \tau}a_{i}v_{i}, \mbox{all}~~ a_{i}\neq
0$. This follows if we take $g_{1}=\Sigma_{v_{i}\subset
\tau_{1}}a_{i}v_{i} (a_{i}\neq 0)$ and $g_{2}=\Sigma_{v_{i}\subset
\tau_{2}}a_{i}v_{i} (a_{i}\neq 0)$. From theorem 4.3.2, $X_{({\bf
g})}=\overline{O}_{\tau}=X_{(g)}$.

If $\tau_{1}$ and  $\tau_{2}$ do not satisfy  the above two types
of conditions, without loss of generality, we suppose
$$\tau_{1}=(v_{0},\cdots,v_{s}), \tau_{2}=(v_{0},\cdots,v_{j},v_{s+1},\cdots,v_{t})
, j<s, t>s.$$ then let $\tau=\tau_{1}[1]\cup
\tau_{2}[1]=(v_{0},\cdots,v_{j},\cdots,v_{s},v_{s+1},\cdots,v_{t})$.
From  proposition 4.3.3, we know $\overline{O}_{\tau}={\bf
P}(Q_{\tau})$, where $Q_{\tau}=(0,\cdots,0,q_{t+1},\cdots,q_{n})$.
While $\overline{O}_{\tau_{1}}={\bf P}(Q_{\tau_{1}})$, where
 $Q_{\tau_{1}}=(0,\cdots,0,q_{s+1},\cdots,q_{t},q_{t+1},\cdots,q_{n})$,
$\overline{O}_{\tau_{2}}={\bf P}(Q_{\tau_{2}})$, where
$Q_{\tau_{2}}=(0,\cdots,0,q_{j+1},\cdots,q_{s},0,\cdots,
0,q_{t+1},\cdots,q_{n})$. Let
$$d_{1}=gcd(q_{s+1,\cdots,q_{n}}),~~~ d_{2}=gcd(q_{j+1},\cdots,q_{s},q_{t+1},\cdots,q_{n})
  \eqno{(5.1)}$$
So from theorem 4.3.5, $(q_{s+1},\cdots,q_{n})$ and
$(q_{j+1},\cdots,q_{s},q_{t+1},\cdots,q_{n})$ are the maximal
subsets of $(q_{0},\cdots,q_{n})$ that satisfy the condition
(5.1). We conclude that $gcd(q_{t+1},\cdots,q_{n})\geq d_{1}d_{2}$
and that $(q_{t+1},\cdots,q_{n})$ must be the maximal subset of
$(q_{0},\cdots,q_{n})$ that satisfies this condition. So from
theorem 4.3.5., $\overline{O}_{\tau}$ is a twisted sector.
\end{proof}

\begin{rem} From the above theorem, every
3-multisector of weighted projective space ${\bf P}(Q)$ of type
$Q=(q_{0},\cdots,q_{n})$ is actually a twisted sector, so we can
refer to the Remark 4.3.6. to describe the orbifold structure of
the 3-multisectors.
\end{rem}
\subsection{The Chen-Ruan Cohomology Rings of
Weighted Projective Spaces.}

In this section we discuss the key point of computing the ring
structure of Chen-Ruan cohomology of weighted projective space
${\bf P}(Q)$ of type $Q=(q_{0},\cdots,q_{n})$, i.e., the orbifold
cup product. The most important part for the orbifold cup product
is the obstruction bundle which was constructed as follows.

Let $X_{({\bf g})}$ be a 3-multisector of $X={\bf P}(Q)$, ${\bf
g}=(g_{1},g_{2},g_{3})\in T_{3}^{0}$. Let $E_{({\bf
g})}\longrightarrow X_{({\bf g})}$ be the obstruction bundle
defined in the section 2.3.1. On the local chart $(V_{x}^{{\bf
g}},C({\bf g}),\pi)$ of $X_{({\bf g})}$, $E_{({\bf g})}$ is given
by $(TV_{x}\otimes H^{0,1}(\Sigma))^{K(\bf{g})}\times
V_{x}^{\bf{g}}\longrightarrow V_{x}^{\bf{g}}$, and the rank of
$E_{({\bf g})}$ is given by the formula(1.3).
$$rank_{\bf{C}}(E_{(\bf{g})})=dim_{\bf{C}}(X_{(\bf{g})})-dim_{\bf{C}}(X)+
\Sigma_{j=1}^{3}\iota_{(g_{j})}
\eqno{(5.2)}$$ \\
The orbifold cup product is defined by: let $\eta_{j}\in
H^{d_{j}}(X_{(g_{j})};\mbox{\bf{Q}})$, for $j=1,2,3$. Define maps
$e_{j}: X_{(\bf{g})}\longrightarrow X_{(g_{j})}$ by
$(x,(\mbox{\bf{g}})_{x})\longmapsto (p,(g_{j})_{x})$. Then

$$<\eta_{1},\eta_{2},\eta_{3}>_{orb}=\int_{X_{(\bf{g})}}^{orb}e^{*}_{1}\eta_{1}\wedge
 e^{*}_{2}\eta_{2}\wedge e^{*}_{3}\eta_{3}\wedge e_{A}(E_{(\bf{g})})
 \eqno{(5.3)}$$
Where $e_{A}(E_{(\bf{g})})$ is the Euler form computed from the
connection $A$.
$$<\eta_{1}\cup_{orb}\eta_{2},\eta_{3}>_{orb}=<\eta_{1},\eta_{2},\eta_{3}>_{orb}
 \eqno{(5.4)}$$
From the above, we have the following proposition.

\begin{prop} Let $\alpha\in
H^{*}_{orb}(X_{(g_{1})};{\bf Q})$, $\beta\in
H^{*}_{orb}(X_{(g_{2})};{\bf Q})$, ${\bf g}=(g_{1},g_{2},g_{3})\in
T_{3}^{0}$, if $\Sigma_{j=1}^{3}\iota_{(g_{j})}>n$, then
$\alpha\cup_{orb}\beta=0$.
\end{prop}

\begin{proof} From (5.2), we have
$\Sigma_{j=1}^{3}\iota_{(g_{j})}-n=rank_{\bf{C}}(E_{(\bf{g})})-dim_{\bf{C}}(X_{(\bf{g})})$.
If $\Sigma_{j=1}^{3}\iota_{(g_{j})}>n$, then
$rank_{\bf{C}}(E_{(\bf{g})})>dim_{\bf{C}}(X_{(\bf{g})})$, so the
integration (5.3) is zero,   $\alpha\cup_{orb}\beta=0$.
\end{proof}

Now in the next three sections we concretely discuss how to
compute the 3-point function defined in (5.3).

\subsection{$\bf{q_{0},\cdots,q_{n}}$ are mutually
prime-A Simple Case.}  Given a $n+1$-tuple
$Q=(q_{0},\cdots,q_{n})$. $q_{i}'s$ are mutually prime. Let ${\bf
P}(Q)$ be the weighted projective space of type $Q$, then the
orbifold singularities are the $n+1$ isolated points:
$p_{i}=[0,\cdots,i,\cdots 0] (i=0,1,\cdots,n)$ with local orbifold
groups ${\bf Z}_{q_{i}} (i=0,1,\cdots,n)$. If we let
$c_{0},\cdots,c_{n}$ be the generators of ${\bf
Z}_{q_{0}},\cdots,{\bf Z}_{q_{n}}$ respectively, then we have
$q_{0}-1$ twisted sectors isomorphic to
$X_{(c_{0})}=p_{0}$,$\cdots$,$q_{n}-1$ twisted sectors isomorphic
to $X_{(c_{n})}=p_{n}$. And we can also see that the 3-sectors are
all isolated points.

If we have $\alpha\in H^{*}(X_{(g_{1})};{\bf Q})$,$\beta\in
H^{*}(X_{(g_{2})};{\bf Q})$, then
$X_{(g_{1},g_{2},(g_{1}g_{2})^{-1})}=\{pt\}$ iff $g_{1},g_{2}$
belong to some ${\bf Z}_{q_{i}} (i=0,1,\cdots,n)$. Without loss of
generality, we assume $g_{1},g_{2}\in {\bf Z}_{0}$, then from the
formula (4.1.7) in [6],
$$\alpha\cup_{orb}\beta=\sum_{(h_{1},h_{2}),~\\ h_{i}\in(g_{i})}
(\alpha\cup_{orb}\beta)_{(h_{1},h_{2})} \eqno{(5.5)}$$ where
$$<(\alpha\cup_{orb}\beta)_{(h_{1},h_{2})},\gamma>_{orb}=
\int_{X_{(h_{1},h_{2})}}e_{1}^{*}\alpha\wedge e_{2}^{*}\beta\wedge
e_{3}^{*}\gamma \wedge e(E_{({\bf g})})  \eqno{(5.6)}
$$
$e_{i}: X_{({\bf g})}\longrightarrow X_{(g_{i})}$ is the map
mentioned above and $E_{({\bf g})}$ is the obstruction bundle over
$X_{({\bf g})}$. From the formula (4.2), the dimension of the
bundle $E_{({\bf g})}$ is:
$$dim(e(E_{({\bf g})}))=2(\iota_{(g_{1})}+\iota_{(g_{2})}+\iota_{(g_{3})})-2n.$$
Because $X_{({\bf g})}$ is a point, so the integration (5.6) is
nonzero iff $\alpha\in H^{0}(X_{(g_{1})};{\bf Q})$,$\beta\in
H^{0}(X_{(g_{2})};{\bf Q})$,$\gamma\in H^{0}(X_{(g_{3})};{\bf Q})$
and $dim(e(E_{({\bf g})}))=0$. At this moment
$\iota_{(g_{1})}+\iota_{(g_{2})}+\iota_{(g_{3})}=n$. Suppose
$\alpha$ and $\beta$ are all generators. Let $\gamma\in
H^{0}(X_{(g_{3})};{\bf Q})$ is the generator. So the integration
(5.6) is:
$$<(\alpha\cup_{orb}\beta)_{(g_{1},g_{2})},\gamma>_{orb}=
\frac{1}{|{\bf Z}_{q_{0}}|}\int_{\{pt\}}e_{1}^{*}\alpha\wedge
e_{2}^{*}\beta\wedge e_{3}^{*}\gamma =\frac{1}{q_{0}}$$ If we let
the $\delta$ be the generator of $H^{0}(X_{(g_{1}g_{2})};{\bf
Q})$, then the integration
$$<\delta,\gamma>_{orb}=\int_{X_{(g_{1}g_{2})}}\delta\wedge
 I^*{\gamma}=\frac{1}{q_{0}}$$
So we have $$\alpha\cup_{orb}\beta=\delta. \eqno{(5.7)}$$

\begin{exa} For $Q=(2,3,5)$, ${\bf P}(Q)={\bf
P}_{2,3,5}^{2}$, we have $q_{0}=2,q_{1}=3,q_{2}=5$. From
proposition 3.2.2., we have
$$
C_{0}= \left(
\begin{array}{cc}
2&1\\
0&1
\end{array}
\right)
$$
So let $v_{1}=(2,1), v_{2}=(0,1)$, we have
$v_{0}=-\frac{q_{1}}{q_{0}}v_{1}-\frac{q_{2}}{q_{0}}v_{2}=(-3,-4)$.
The fan of ${\bf P}_{2,3,4}^{2}$ is generated by
$\{v_{0},v_{1},v_{2}\}$. For
$\sigma_{2}=(v_{0},v_{1})=((-3,-4),(2,1))$, we have
$G_{\sigma_{2}}=N/N_{\sigma_{2}}={\bf Z}_{5}$. We  write the
matrix representation of the action of ${\bf Z}_{5}$ on
$U_{\sigma_{2}^{'}}={\bf C}^{2}$ as follows.
$$
\left(
\begin{array}{cc}
1&0\\
0&1
\end{array}
\right); \left(
\begin{array}{cc}
e^{2\pi i\cdot \frac{1}{5}}&0\\
0&e^{2\pi i\cdot \frac{4}{5}}
\end{array}
\right) ; \left(
\begin{array}{cc}
e^{2\pi i\cdot \frac{2}{5}}&0\\
0&e^{2\pi i\cdot \frac{3}{5}}
\end{array}
\right)
$$
$$
\left(
\begin{array}{cc}
e^{2\pi i\cdot \frac{3}{5}}&0\\
0&e^{2\pi i\cdot \frac{2}{5}}
\end{array}
\right); \left(
\begin{array}{cc}
e^{2\pi i\cdot \frac{4}{5}}&0\\
0&e^{2\pi i\cdot \frac{1}{5}}
\end{array}
\right)
$$
For $\sigma_{1}=((-3,-4),(0,1))$, we have $G_{\sigma_{1}}={\bf
Z}_{3}$. We  write the matrix representation of the action of
${\bf Z}_{3}$ on $U_{\sigma_{1}^{'}}={\bf C}^{2}$ as follows.
$$
\left(
\begin{array}{cc}
1&0\\
0&1
\end{array}
\right); \left(
\begin{array}{cc}
e^{2\pi i\cdot \frac{1}{3}}&0\\
0&e^{2\pi i\cdot \frac{1}{3}}
\end{array}
\right) ;  \left(
\begin{array}{cc}
e^{2\pi i\cdot \frac{2}{3}}&0\\
0&e^{2\pi i\cdot \frac{2}{3}}
\end{array}
\right)
$$
And for $\sigma_{0}=(v_{1},v_{2})=((2,1),(0,1))$, we have
$G_{\sigma_{0}}={\bf Z}_{2}$ and representation:
$$
\left(
\begin{array}{cc}
1&0\\
0&1
\end{array}
\right);  \left(
\begin{array}{cc}
e^{2\pi i\cdot \frac{1}{2}}&0\\
0&e^{2\pi i\cdot \frac{1}{2}}
\end{array}
\right)
$$
If we let
$$ g_{1}=\left(
\begin{array}{cc}
e^{2\pi i\cdot \frac{1}{5}}&0\\
0&e^{2\pi i\cdot \frac{4}{5}}
\end{array}
\right); g_{2}=\left(
\begin{array}{cc}
e^{2\pi i\cdot \frac{1}{3}}&0\\
0&e^{2\pi i\cdot \frac{1}{3}}
\end{array}
\right); g_{3}=\left(
\begin{array}{cc}
e^{2\pi i\cdot \frac{1}{2}}&0\\
0&e^{2\pi i\cdot \frac{1}{2}}
\end{array}
\right)
$$
then we have the twisted  sectors: let
$p_{0}=[1,0,0],p_{1}=[0,1,0], p_{2}=[0,0,1]$,
$X_{(g_{1})}=X_{(g_{1}^{2})}=X_{(g_{1}^{3})}=X_{(g_{1}^{4})}=p_{0}$,
$X_{(g_{2})}=X_{(g_{2}^{2})}=p_{1}$ and $X_{(g_{3})}=p_{3}$.\\

Since  $$
\begin{array}{l}
\iota_{(g_{1})}+\iota_{(g_{1}^{4})}=2,\\
\iota_{(g_{1}^{2})}+\iota_{(g_{1}^{3})}=2,\\
\iota_{(g_{2})}+\iota_{(g_{2})}+\iota_{(g_{2})}=2,\\
\iota_{(g_{2})}+\iota_{(g_{2}^{2})}=2,\\
\iota_{(g_{3})}+\iota_{(g_{3})}=2.
\end{array}
$$

Let $\alpha_{1},\alpha_{2},\alpha_{3},\alpha_{4}$ be the
generators of $H^{0}(X_{(g_{1})};{\bf
Q}),H^{0}(X_{(g_{1}^{2})};{\bf Q}),  H^{0}(X_{(g_{1}^{3})};{\bf
Q}),$ $H^{0}(X_{(g_{1}^{4})};{\bf Q})$ respectively,
$\beta_{1},\beta_{2}$ be the generators of $H^{0}(X_{(g_{2})};{\bf
Q}), H^{0}(X_{(g_{2}^{2})};{\bf Q})$, and $\gamma$ be the
generators of $H^{0}(X_{(g_{3})};{\bf Q})$. We also let $e_{0}$ be
the generator of $H^{0}({\bf P}^{2}_{2,3,5};{\bf Q})$. So from the
above discussion in section 5.4. and the formula (5.7), we have
$$
\begin{array}{l}
\alpha_{1}\cup_{orb}\alpha_{4}=e_{0},\\
\alpha_{2}\cup_{orb}\alpha_{3}=e_{0},\\
\beta_{1}\cup_{orb}\beta_{1}=\beta_{2},\\
\beta_{1}\cup_{orb}\beta_{2}=e_{0},\\
\gamma\cup_{orb}\gamma=e_{0}.$$
\end{array}
$$\\
\end{exa}

\subsection{The Obstruction Bundle.}  In this section we
introduce a method to compute the obstruction bundle locally. Let
$E_{({\bf g}}\longrightarrow X_{({\bf g})}$ be the obstruction
bundle over the 3-multisector $X_{({\bf g})}$. For a weighted
projective space ${\bf P}(Q)$ of type $Q=(q_{0},\cdots,q_{n})$.
From theorem 5.2.2, every 3-multisector $X_{({\bf g})}$ is a
twisted sector.  Assume $X_{({\bf g})}=\overline{O}_{\tau}={\bf
P}(Q_{\tau})$, $\tau=(v_{0},\cdots,v_{i-1})$ is a cone in the fan
$\Xi$, $Q_{\tau}=(q_{i},\cdots,q_{n})$ and
$gcd(q_{i},\cdots,q_{n})=d>1$. So we have
$g_{1},g_{2},g_{3}\in{\bf Z}_{d}$ and $g_{1}g_{2}g_{3}=1$. ${\bf
P}(Q_{\tau})={\bf P}^{n}(0,\cdots,0,q_{i},\cdots,q_{n})$ is a
hyperplane of ${\bf P}(Q)$. If we let
$o(g_{1})=l,o(g_{2})=q,o(g_{3})=r$, then $l+q+r\equiv 0(mod~ d)$.
From section 2.3., the fundamental group of orbifold
$S^{2}(l,q,r)$ is:
$$\pi_{1}^{orb}(S^{2})=\{\lambda_{1},\lambda_{2},\lambda_{3}|\lambda_{1}^{l}=
\lambda_{2}^{q}=\lambda_{3}^{r}=\lambda_{1}\lambda_{2}\lambda_{3}=1\}$$
Then $\Gamma=\pi_{1}^{orb}(S^{2})$ is a fuchsian group. Let
$\varphi: \Gamma\longrightarrow {\bf Z}_{d}$ be the homomorphism
defined by $\lambda_{i}\longrightarrow g_{i}$ for $i=1,2,3$.
Assume $Im\varphi=G=K({\bf g})$,  there is a Riemann surface
$\Sigma$ such that $\Sigma/K({\bf g})=S^{2}(l,q,r)$. From [16],
since the group $K({\bf g})$ is also a cyclic group, we can
describe the $(l,q,r)$ as three cases:
$$
\begin{array}{ll}
  (1)&(l,q,r)=(p,p,0), p>1;\\
  (2)&(l,q,r)=(3,3,3), (2,3,6),(2,4,4);\\
  (3)&(l,q,r)~~ \mbox{is not the cases of (1) and (2)}.
\end{array}
$$
From (1) we can see that the Riemann surface $\Sigma=S^{2}$ and
the cyclic group ${\bf Z}_{p}$ acts on $S^{2}$ with quotient
$S^{2}(p,p,0)$. In the case (2), the three groups of values are
the only choices such that the Riemann surface $\Sigma$ is the
torus. In this case, $\frac{1}{p}+\frac{1}{q}+\frac{1}{r}=1$. In
the case (3), we know that the universal covering space of
$S^{2}(l,q,r)$ is ${\bf H}^{2}$(Poincare disk model),then
$\Sigma={\bf H}^{2}/ker(\varphi)$. By Riemann-Hurwitz formula,
$$g(\Sigma)=\frac{|K({\bf g})|+2-|K({\bf g})|(\frac{1}{p}+\frac{1}{q}+\frac{1}{r})}
{2} \eqno{(5.8)}$$ In this case
$\frac{1}{p}+\frac{1}{q}+\frac{1}{r}>1$ and $g(\Sigma)\geq 2$.

Now we assume that the genus of $\Sigma$ is $g$. Let
$\gamma_{1},\cdots,\gamma_{2g}$ is a canonical homological basis
satisfying the following relations:
$$
\left\{\begin{array}{ll}
\gamma_{i}\gamma_{g+i}=1, \forall ~~i=1,\cdots,g; \\
\gamma_{i}\gamma_{j}=0, ~~~~\mbox{otherwise}.
\end{array} \right.
$$
Then by Riemann bilinear relation [15], the period matrix is:
${I_{g} \choose Z}$ such that:

$$
\left\{
\begin{array}{ll}
a.)& Z=^{t}Z;\\
b.)& ImZ ~~\mbox{is positively definitive}.
\end{array}
\right.
$$

Suppose $\check{\gamma}_{1},\cdots,\check{\gamma}_{2g}$ is the
dual basis in $H^{1}(\Sigma)$, then there exists a basis
$\omega_{1},\cdots,\omega_{g}$ in $H^{0,1}(\Sigma)$ such that
$$(\omega_{1},\cdots,\omega_{g})=(\check{\gamma}_{1},\cdots,\check{\gamma}_{2g})
{I_{g} \choose Z}$$

The element $\lambda_{1}$ in $\Gamma$ acts on $H^{1}(\Sigma)$
naturally. Let $A^{t}$ be the representation matrix of
$(\lambda_{1})_{*}$ under the basis
$\gamma_{1},\cdots,\gamma_{2g}$, then $\lambda_{1}^{*}$ acts on
$H^{1}(\Sigma)$, the matrix representation is $A$. So:
$$\lambda_{1}^{*}(\omega_{1},\cdots,\omega_{g})=
(\check{\gamma}_{1},\cdots,\check{\gamma}_{2g})\cdot A\cdot {I_{g}
\choose Z}.$$ We know that $A\cdot {I_{g} \choose Z}$ is a
$2g\times g$ matrix. Assume the matrix which consists of the first
$g$ rows of the matrix $A\cdot {I_{g} \choose Z}$ is $U$. Since
$\lambda_{1}$ is a holomorphic map, $\lambda_{1}^{*}$ preserves
the subspace $H^{1,0}(\Sigma)$. So
$$\lambda_{1}^{*}(\omega_{1},\cdots,\omega_{g})=
(\omega_{1},\cdots,\omega_{g})\cdot U$$ We have
$$\lambda_{1}^{*}(\omega_{1},\cdots,\omega_{g})=
(\check{\gamma}_{1},\cdots,\check{\gamma}_{2g})\cdot {I_{g}
\choose Z}\cdot U$$ Thus we can determine the matrix $Z$, and the
matrix $U$ can also be determined. We know that
$H^{0,1}(\Sigma)=\overline{H^{1,0}(\Sigma)}$, so the basis of
$\overline{H^{1,0}(\Sigma)}$ is
$\overline{\omega}_{1},\cdots,\overline{\omega}_{g}$, and
$$\lambda_{1}^{*}(\overline{\omega}_{1},\cdots,\overline{\omega}_{g})=
(\overline{\omega}_{1},\cdots,\overline{\omega}_{g})\cdot
\overline{U} \eqno{(5.9)}$$

Let $U_{j}=\{[z]_{Q}\in {\bf P}^{n}_{q_{0},\cdots,q_{n}}:
z_{j}\neq 0\}\subset {\bf P}^{n}_{q_{0},\cdots,q_{n}}$ for
$j=0,\cdots,n$. Then since $Q_{\tau}=(q_{i},\cdots,q_{n})$, we see
that $X_{({\bf g})}={\bf P}(Q_{\tau})$ can be covered by $X_{({\bf
g})}\cap U_{i},\cdots,X_{({\bf g})}\cap U_{n}$. From section 3.1,
for $j\geq i$, we have a bijective map $\phi_{j}$ from $U_{j}$ to
$\mbox{\bf{C}}^{n}/\mu_{q_{j}}(Q_{q_{j}})$ given by
$$\phi_{j}([z]_{Q})=\left(\frac{z_{0}}{(z_{j})^{q_{0}/q_{j}}},\cdots,
\frac{\hat{z_{j}}}{z_{j}},\cdots,\frac{z_{n}}{(z_{j})^{q_{n}/q_{j}}}\right)_{q_{j}}$$
So we choose the coordinates of ${\bf C}^{n}=V_{j}$ by
$\left(\frac{z_{0}}{(z_{j})^{q_{0}/q_{j}}},\cdots,
1,\cdots,\frac{z_{n}}{(z_{j})^{q_{n}/q_{j}}}\right)$. If we let
$x_{0}=\frac{z_{0}}{(z_{j})^{q_{0}/q_{j}}},\cdots,x_{j}=1,
\cdots,x_{n}=\frac{z_{n}}{(z_{j})^{q_{n}/q_{j}}}$, then let
$p_{j}=[0,\cdots,1,\cdots,0]$ be the point in $V_{j}$,
$(TV_{j})_{p_{j}}$ has the basis $\frac{\partial}{\partial
x_{0}},\cdots, \frac{\partial}{\partial x_{n}}$, and the
$\lambda_{1}^{*}=g_{1}$ acts on $(TV_{j})_{p_{j}}$ in the natural
way. We denote the diagonal representation matrix of
$\lambda_{1}^{*}$ by $D$. So on $(H^{0,1}(\Sigma)\otimes
(TV_{j})_{p_{j}})$ we have a basis:
$$\left\{\frac{\partial}{\partial
x_{0}}\otimes\overline{\omega}_{1},\frac{\partial}{\partial
x_{1}}\otimes\overline{\omega}_{1},\cdots,\frac{\partial}{\partial
x_{n}}\otimes\overline{\omega}_{1},\frac{\partial}{\partial
x_{0}}\otimes\overline{\omega}_{2},\cdots,\frac{\partial}{\partial
x_{n}}\otimes\overline{\omega}_{n}\right\}.$$ and
$\lambda_{1}^{*}=g_{1}$ acts on $(H^{0,1}(\Sigma)\otimes
(TV_{j})_{p_{j}})$, which the matrix representation is $D\otimes
\overline{U}$. Because $g_{1},g_{2},g_{3}$ generate the cyclic
group ${\bf Z}_{d}$, and assume that the $g_{1}$ is a generator,
so the matrix $D\otimes \overline{U}$ has eigenvalue $1$ with
multiplicity $e$, where
$$e=dim_{{\bf C}}(e(E_{({\bf g})}))=dim_{{\bf C}}(X_{({\bf g})})
-n+\Sigma_{j=1}^{3}\iota_{(g_{j})} \eqno{(5.10)}$$ So locally the
obstruction bundle is generated by the $e$ eigenvectors
$\xi_{1},\cdots,\xi_{e}$.

Suppose $\xi=\alpha\otimes \omega\in \left((TV_{j})_{p_{j}}\otimes
H^{0,1}(\Sigma)\right)^{K({\bf g})}$, then $\alpha$ is a linear
combination of $\frac{\partial}{\partial
x_{0}},\cdots,\frac{\partial}{\partial x_{i-1}}$, otherwise
$\omega$ will be a harmonic form on $S^{2}$ because $\omega$ is
$K({\bf g})$-invariant, so $\omega=0$. And also from the action of
$\lambda_{1}^{*}=g_{1}$ on $(TV_{j})_{p_{j}}\otimes
H^{0,1}(\Sigma)$ and $C({\bf g})$ on the invariant subspace
$((TV_{j})_{p_{j}}\otimes H^{0,1}(\Sigma))^{{\bf Z}_{d}}$, we  see
that each eigenvector $\xi_{t}$ can be written as the linear
combination of $\{\frac{\partial}{\partial x_{t}}\otimes
\overline{\omega}_{1},\cdots,\frac{\partial}{\partial
x_{t}}\otimes \overline{\omega}_{g}\}$ for some $t$ with $0\leq
t\leq i-1$.

Now assume $k\geq i$, and $U_{k}\cap X_{({\bf g})}$ is another
open subset of $X_{({\bf g})}$. Let $U_{k}=V_{k}/{\bf Z}_{k}$,
from the above discussion, we can choose the coordinates of
$V_{k}$ as:
$$\left(y_{0}=\frac{z_{0}}{(z_{k})^{q_{0}/q_{k}}},\cdots,
y_{k}=1,\cdots,y_{n}=\frac{z_{n}}{(z_{k})^{q_{n}/q_{k}}}\right).$$
So we have a basis on $(H^{0,1}(\Sigma)\otimes (TV_{k})_{p_{k}})$:

$$\left\{\frac{\partial}{\partial
y_{0}}\otimes\overline{\omega}_{1},\frac{\partial}{\partial
y_{1}}\otimes\overline{\omega}_{1},\cdots,\frac{\partial}{\partial
y_{n}}\otimes\overline{\omega}_{1},\frac{\partial}{\partial
y_{0}}\otimes\overline{\omega}_{2},\cdots,\frac{\partial}{\partial
y_{n}}\otimes\overline{\omega}_{n}\right\}.$$

Then we see that on the neighborhood $U_{k}\cap X_{({\bf g})}$,
the obstruction bundle is still generated by $e$ eigenvectors
$\xi_{1}^{'},\cdots,\xi_{e}^{'}$. And $\xi_{l}^{'}$ can be
obtained from $\xi_{l}$ by substituting $\frac{\partial}{\partial
y_{l}}$ for $\frac{\partial}{\partial x_{l}}$. On the $V_{k}$, we
have:

\begin{eqnarray}
\frac{\partial}{\partial x_{0}} =\Sigma\frac{\partial
y_{j}}{\partial x_{0}} \frac{\partial}{\partial y_{j}}
&=&\frac{\partial}{\partial
x_{0}}(\frac{x_{0}}{(x_{k})^{q_{0}/q_{k}}})
\frac{\partial}{\partial y_{0}} \nonumber \\
&=&\frac{1}{(x_{k})^{q_{0}/q_{k}}}\frac{\partial}{\partial y_{0}} \nonumber \\
&=&(y_{j})^{q_{0}/q_{j}}\frac{\partial}{\partial y_{0}} \nonumber
\end{eqnarray}

For the same computation, we also have:
 $$
\left\{\begin{array}{l}
\frac{\partial}{\partial x_{1}}=(y_{j})^{q_{1}/q_{j}}\frac{\partial}{\partial y_{1}};\\
\cdots \cdots \\
\frac{\partial}{\partial
x_{i-1}}=(y_{j})^{q_{i-1}/q_{j}}\frac{\partial}{\partial y_{i-1}}.
\end{array}
\right.
$$

Assume on the neighborhood $U_{j}\cap X_{({\bf g})}$, the
eigenvectors $\xi_{1},\cdots,\xi_{e}$ are generated by the linear
combination of the tuples $\{\frac{\partial}{\partial
x_{t_{1}}}\otimes
\overline{\omega}_{1},\cdots,\frac{\partial}{\partial
x_{t_{1}}}\otimes
\overline{\omega}_{g}\},\cdots,\{\frac{\partial}{\partial
x_{t_{e}}}\otimes
\overline{\omega}_{1},\cdots,\frac{\partial}{\partial
x_{t_{e}}}\otimes \overline{\omega}_{g}\}$ respectively for $0\leq
t_{1},\cdots,t_{e}\leq i-1$. Then  the transition function:
$$h_{kj}: (U_{j}\cap X_{({\bf g})})\times {\bf C}^{e}\hookleftarrow
(U_{j}\cap U_{k}\cap X_{({\bf g})})\times {\bf
C}^{e}\longrightarrow (U_{j}\cap U_{k}\cap X_{({\bf g})})\times
{\bf C}^{e}\hookrightarrow (U_{k}\cap X_{({\bf g})})\times {\bf
C}^{e}$$ can be written as:
$$h_{kj}(x;c_{1},\cdots,c_{e})=\left(x;y_{j}^{q_{t_{1}}/q_{j}}(x)\cdot
c_{1},\cdots, y_{j}^{q_{t_{e}}/q_{j}}(x)\cdot c_{e}\right).
\eqno{(5.11)}$$ So locally the transition matrix is:
 $$
\left(
\begin{array}{ccc}
y_{j}^{q_{t_{1}}/q_{j}}(x)&&\\
&\ddots&\\
&&y_{j}^{q_{t_{e}}/q_{j}}(x)
\end{array}
\right)
$$
Then the obstruction bundle $E_{({\bf g})}$ can be splited as the
whitney sum of line bundles, let $E_{({\bf
g})}=\oplus_{l=1}^{e}E_{l}$. Every line bundle $E_{l}$ is
generated by $\xi_{l}$ on the neighborhood $U_{j}\cap X_{({\bf
g})}$. The group ${\bf Z}_{d}$ acts diagonally on the obstruction
bundle $E_{({\bf g})}$, so it acts on every line bundle $E_{l}$
naturally. Assume the matrix representation of the action of the
generator of ${\bf Z}_{d}$ on the obstruction bundle $E_{({\bf
g})}$ is:
$$
\left(
\begin{array}{ccc}
e^{2\pi i\cdot \frac{m_{1}}{d}}&&\\
&\ddots&\\
&&e^{2\pi i\cdot \frac{m_{e}}{d}}
\end{array}
\right)
$$
here $0\leq m_{l}< d, 1\leq l\leq e$. Then we have the following
facts. $e^{2\pi i\cdot \frac{m_{l}}{d}}$ is a $d_{l}$-root of $1$
for $1\leq l\leq e$, and it is clear that $d_{l}$ is a divisor of
$d$.

\subsection{Computation of the 3-Point Function.}  In this
section we use the localization technique ([3],[12]) to calculate
the 3-point function defined in the orbifold cup product.

Let $X={\bf P}(Q)$ be the weighted projective space of type
$Q=(q_{0},\cdots,q_{n})$ and $X_{({\bf g})}$ be a 3-multisector.
Then $X_{({\bf g})}$ is a twisted sector from theorem 5.2.2.,
assume $X_{({\bf g})}=\overline{O}_{\tau}={\bf P}(Q_{\tau})$,
where $Q_{\tau}=(0,\cdots,0,q_{i},\cdots,q_{n})$,${\bf
g}=(g_{1},g_{2},g_{3})\in T_{3}^{0}$ and
$\tau=(v_{0},\cdots,v_{i-1})$ is a cone of the fan
$\Xi=\{v_{0},\cdots,v_{n}\}$. The orbifold structure of $X_{({\bf
g})}$ can be described  in Remark 4.3.6., From (2.4), the key
calculation of the orbifold cup product is to calculate the
3-point function:
$$<\eta_{1},\eta_{2},\eta_{3}>_{orb}=\int_{X_{(\bf{g})}}^{orb}e^{*}_{1}\eta_{1}\wedge
 e^{*}_{2}\eta_{2}\wedge e^{*}_{3}\eta_{3}\wedge e_{A}(E_{(\bf{g})})
 \eqno{(5.12)}$$
where $\eta_{j}\in H^{*}(X_{(g_{j})};\mbox{\bf{Q}})$, for
$j=1,2,3$.

Now we analyze the integral formula (5.12). In order to compute
conveniently, we prove that we always can suppose $q_{0}=1$. If
$q_{0}\neq 1$, let $\widetilde{Q}=(1,q_{0},\cdots,q_{n})$, then
${\bf P}^{n}_{q_{0},\cdots,q_{n}}\subset {\bf
P}(\widetilde{Q})={\bf P}^{n+1}_{1,q_{0},\cdots,q_{n}}=Y$ is a
hypersurface which is obtained if by letting the first homogeneous
coordinate of ${\bf P}(\widetilde{Q})$ be zero. From theorem 4.3.5
and theorem 5.2.2, it is easy to see that ${\bf P}(Q)$ and  ${\bf
P}(\widetilde{Q})$ have the same twisted sectors and
3-multisectors. Suppose that the matrix representations in ${\bf P
}(\widetilde{Q})$ corresponding to $g_{1},g_{2},g_{3}$ in ${\bf P
}(Q)$ are $\widetilde{g}_{1},\widetilde{g}_{2},\widetilde{g}_{3}$,
then $Y_{(\widetilde{{\bf g
}})}=Y_{(\widetilde{g}_{1},\widetilde{g}_{2},\widetilde{g}_{3})}=X_{({\bf
g })}$. The cohomological classes
$e_{1}^{*}\eta_{1},e_{2}^{*}\eta_{2},e_{3}^{*}\eta_{3}$ are
invariant when they are taken as the cohomological classes of
$Y_{(\widetilde{{\bf g}})}$. Suppose the homogeneous coordinates
of ${\bf P}(\widetilde{Q})$ is ${\bf z}=[z,z_{0},\cdots,z_{n}]$,
let $\widetilde{U}_{j}=\{z_{j}\neq 0|{\bf z}\in {\bf
P}(\widetilde{Q})\}$, $(0\leq j\leq n)$, $\widetilde{U}=\{z\neq
0|{\bf z}\in {\bf P}(\widetilde{Q})\}$, then $Y_{(\widetilde{{\bf
g}})}$ can be covered by $\bigsqcup_{j=i}^{n}\widetilde{U}_{j}\cap
Y_{(\widetilde{{\bf g}})}$. For the local chart
$\widetilde{U}_{j}$, let $\widetilde{U}_{j}=\widetilde{V}_{j}/{\bf
Z}_{q_{j}}$ and choose the coordinates of $\widetilde{V}_{j}$ as
$(x=\frac{z}{(z_{j})^{1/q_{j}}},x_{0}=\frac{z_{0}}{(z_{j})^{q_{0}/q_{j}}},\cdots,
x_{n}=\frac{z_{n}}{(z_{j})^{q_{n}/q_{j}}})$, so we have a base of
$(T\widetilde{V}_{j})_{p_{j}}$:  $(\frac{\partial}{\partial
x},\frac{\partial}{\partial x_{0}},\cdots,\frac{\partial}{\partial
x_{n}})$. Because the invariant subspace $((TV_{j})_{p_{j}}\otimes
H^{0,1}(\Sigma))^{{\bf Z}_{d}}$ is generated by
$\xi_{1},\cdots,\xi_{e}$, and we can see that the space
$((TV_{j})_{p_{j}}\otimes H^{0,1}(\Sigma))^{{\bf Z}_{d}} \subset
((T\widetilde{V}_{j})_{\widetilde{p}_{j}}\otimes
H^{0,1}(\Sigma))^{{\bf Z}_{d}}$. We construct a new obstruction
bundle $E_{(\widetilde{{\bf g}})}$ over $Y_{(\widetilde{{\bf
g}})}$ as follows. On the local chart $\widetilde{U}_{j}\cap
Y_{(\widetilde{{\bf g}})}$, this bundle is given by
$\widetilde{V}_{j}\cap H\times ((TV_{j})_{p_{j}}\otimes
H^{0,1}(\Sigma))^{{\bf Z}_{d}}\longrightarrow
\widetilde{V}_{j}\cap H$, where $H=\{x=x_{0}=\cdots=x_{i-1}=0|{\bf
x}\in \widetilde{V}_{j}\}$ is a hypersurface of
$\widetilde{V}_{j}$.  It is easy to see that the transition
function of this bundle is also given by (5.11), so it can also be
splitted into the Whitney sum of line bundles. It is clear that
$E_{(\widetilde{{\bf g}})}\cong E_{({\bf g})}$, so we have
$$\int_{X_{(\bf{g})}}^{orb}e^{*}_{1}\eta_{1}\wedge
 e^{*}_{2}\eta_{2}\wedge e^{*}_{3}\eta_{3}\wedge
 e_{A}(E_{(\bf{g})})=\int_{Y_{(\widetilde{\bf{g}})}}^{orb}e^{*}_{1}\eta_{1}\wedge
 e^{*}_{2}\eta_{2}\wedge e^{*}_{3}\eta_{3}\wedge
 e_{\widetilde{A}}(E_{(\widetilde{\bf{g}})})
 \eqno{(5.13)}$$
So in the following analysis, we assume that $q_{0}=1$, and we
give a formula to compute the integration (5.12).

First from  Remark 4.3.6, we know that if
$d=gcd(q_{i},\cdots,q_{n})\neq 1$, $X_{({\bf g})}$ is a nonreduced
orbifold. From the discussion of section 5.5., let $E_{({\bf g}
)}=\oplus_{l=1}^{e}E_{l}$, then for every line bundle $E_{l}$,
using the same method of Park and Poddar [22], consider the
associated orbifold principal bundle $P_{l}$ of $E_{l}$ such that
$E_{l}=P_{l}\times_{S^{1}}{\bf C}$. We know that there is a global
action of ${\bf Z}_{d_{l}}$ on each fibre $F=S^{1}$. The quotient
$P_{l}/{\bf Z}_{d_{l}}$ is again an orbifold principal bundle over
the orbifold $X_{({\bf g})}$. Let $\pi_{l}: P_{l}\longrightarrow
P_{l}/{\bf Z}_{d_{l}}$ be the quotient map, which extends to an
orbifold bundle map. Choose an orbifold connection $A_{l}$ that is
the pullback $\pi_{l}^{*}(A_{l}^{'})$, where $A_{l}^{'}$ is an
orbifold connection on the associated bundle
$E_{l}^{'}=(P_{l}/{\bf Z}_{d_{l}})\times_{S^{1}}{\bf C}$. The Lie
algebra of $F$ can be identified with ${\bf R}$, then the induced
map on the lie algebra $(\pi_{l})_{*}: {\bf R}\longrightarrow {\bf
R}$ is just given by $a\longmapsto d_{l}a$.

Let $\Omega_{l}$ and $\Omega_{l}^{'}$ be the curvature 2-forms for
$A_{l}$ and $A_{l}^{'}$. By proposition 6.2 of [19],
$(\pi_{l})^{*}(\Omega_{l}^{'})=d_{l}\Omega_{l}$. So

$\int_{X_{(\bf{g})}}^{orb}e^{*}_{1}\eta_{1}\wedge
e^{*}_{2}\eta_{2}\wedge e^{*}_{3}\eta_{3}\wedge
e_{A}(E_{(\bf{g})})=
\int_{X_{(\bf{g})}}^{orb}e^{*}_{1}\eta_{1}\wedge
e^{*}_{2}\eta_{2}\wedge e^{*}_{3}\eta_{3}\wedge \Pi_{l=1}^{e}
e_{A_{l}}(E_{l})$ $$=\frac{1}{\Pi_{l}
d_{l}}\int_{X_{(\bf{g})}}^{orb}e^{*}_{1}\eta_{1}\wedge
e^{*}_{2}\eta_{2}\wedge e^{*}_{3}\eta_{3}\wedge \Pi_{l=1}^{e}
e_{A_{l}^{'}}(E_{l}^{'}) \eqno{(5.14)}$$

Since the action of ${\bf Z}_{d_{l}}$ in any uniformizing system
of $E_{l}^{'}$ is trivial, $E_{l}^{'}$  induces an orbifold bundle
$E_{l}^{''}$ over the reduced orbifold $X_{({\bf g})}^{'}$ which
has  an induced connection $A_{l}^{''}$. The connections
$A_{l}^{'}$ and $A_{l}^{''}$ may be represented by the same 1-form
over $V$ for $(V\times {\bf C},G^{'}/{\bf
Z}_{d_{l}},\widetilde{\pi}^{''}_{1})$ of $E_{l}^{'}$ and
$E_{l}^{''}$ respectively. By Chern-Weil theory, $\Omega_{l}^{'}$
and $\Omega_{l}^{''}$ can therefore be represented by the same
2-form on $V$. We know that
$e_{1}^{*}(\eta_{1}),e_{2}^{*}(\eta_{2})$ and
$e_{3}^{*}(\eta_{3})$ are invariant when taken as the cohomology
classes of $X_{({\bf g})}^{'}$. Since $K({\bf g})$ acts on
$X_{({\bf g})}$ trivially, so from (2.1), we have:
$$\int_{X_{(\bf{g})}}^{orb}e^{*}_{1}\eta_{1}\wedge
e^{*}_{2}\eta_{2}\wedge e^{*}_{3}\eta_{3}\wedge \Pi_{l=1}^{e}
e_{A_{l}^{'}}(E_{l}^{'})= \frac{1}{|{\bf
Z}_{d}|}\int_{X_{(\bf{g})}^{'}}^{orb}e^{*}_{1}\eta_{1}\wedge
e^{*}_{2}\eta_{2}\wedge e^{*}_{3}\eta_{3}\wedge \Pi_{l=1}^{e}
e_{A_{l}^{''}}(E_{l}^{''}) \eqno{(5.15)}$$

Next we mainly discuss the method to calculate the integration in
(5.15). From above $X_{({\bf g})}={\bf P}(Q_{\tau})$,
$Q_{\tau}=(0,\cdots,0,q_{i},\cdots,q_{n})$. So we have $X_{({\bf
g})}^{'}={\bf P}(Q_{\tau}/d)$,
$Q_{\tau}/d=(q_{i}/d,\cdots,q_{n}/d), d=gcd(q_{i},\cdots,q_{n})$.
Let $N_{\tau}$ be the sublattice of $N$ generated by
$\tau=(v_{0},\cdots,v_{i-1})$ and $N(\tau)=N/N_{\tau}$ be the
quotient lattice. The fan of $X_{({\bf g})}^{'}$ is given by the
projection of $\Xi$ to $N(\tau)\otimes {\bf R}$. The dual lattice
of $N(\tau)$ is $M(\tau)=\tau^{\perp}\cap M$. The torus
$T=spec({\bf C} [M(\tau)])=O_{\tau}$. The characters $\chi^{m}$
correspond to rational functions on $X_{({\bf g})}^{'}$ when $m\in
M(\tau)$. Then we can use the localization technique of [3], when
reduced to orbifold, to calculate  the integration (5.15).

We know $X_{({\bf g})}^{'}$ is a toric variety, $T$ acts on
$X_{({\bf g})}^{'}$  and this action has $n-i+1$ fixed points
$p_{j}$ for $i\leq j\leq n$. We let $\{\rho_{i}\}$ be the basic
characters of $T$ action and $\{\lambda_{i}\}$ be the parameters
of the lie algebra $t_{{\bf C}}$ of $T$ corresponding to the above
base $\{\rho_{i}\}$. Because $q_{0}=1$, it is clear that the
matrix $C_{0}$ defined in proposition 3.2.2 is the unit matrix. We
compute that the fan of the weighted projective space ${\bf P}(Q)$
is generated by $v_{0},\cdots,v_{n}$, where
$v_{0}=(-q_{1},\cdots,-q_{n})$, $v_{j}=e_{j}$ for $1\leq j\leq n$.
Let $\{m_{1},\cdots,m_{n}\}$ be the standard basis of $M$, we
calculate the base of $M(\tau)$ as:
$$\left\{\rho_{1}=\frac{q_{n}}{d}m_{i}-\frac{q_{i}}{d}m_{n},
\cdots, \rho_{r}=\frac{q_{n}}{d}m_{j}-\frac{q_{j}}{d}m_{n},\cdots,
\rho_{n-i}=\frac{q_{n}}{d}m_{n-1}-\frac{q_{n-1}}{d}m_{n}\right\}$$
where $j=i+r-1$. We first study the action of $T$ on the normal
bundle of $p_{j}$, i.e., the orbifold tangent space $(TX_{({\bf
g})}^{'})_{p_{j}}$.

Consider the fixed points $p_{j}(i\leq j\leq n-1)$. Denote the
local coordinates on a uniformizing system of $X_{({\bf g})}^{'}$
around $p_{j}$ by: $[x_{i},\cdots,1,\cdots,x_{n}]$.  Let
$m^{1}=a_{1}\rho_{1}+\cdots+a_{n-i}\rho_{n-i}$, and
$<m^{1},v_{i}>=1$, $<m^{1},v_{k}>=0$ for $k>i,k\neq j$. Then we
have $a_{r}=-\frac{dq_{i}}{q_{n}q_{j}}, a_{1}=\frac{d}{q_{n}}$, so
$\chi^{m^{1}}=x_{i}$. Similarly, we compute $\chi^{m^{t}}=x_{t}$
for $t\neq r,n-i$,
$m^{t}=\frac{d}{q_{n}}\rho_{t}-\frac{dq_{i+t-1}}{q_{n}q_{j}}\rho_{r}$.
Using the same method, let
$m^{n-i}=a_{1}\rho_{1}+\cdots+a_{n-i}\rho_{n-i}$, and
$<m^{n-i},v_{k}>=0$ for $k\geq i,k\neq j,n$, $<m^{n-i},v_{n}>=1$,
we have $a_{r}=-\frac{d}{q_{j}}$,
$m^{n-i}=-\frac{d}{q_{j}}\rho_{r}$, $\chi^{m^{n-i}}=x_{n}$. So the
$T$-equivariant Euler class of normal bundle of $p_{j} (i\leq
j\leq n)$ is given by
$$e_{T}\left(\nu_{p_{j}}\right)=\left(-\frac{d}{q_{j}}\lambda_{r}\right)\prod_{k\neq r}\frac{d}{q_{n}}
\left(\lambda_{k}-\frac{q_{i+k-1}}{q_{j}}\lambda_{r}\right)
  \eqno{(5.16)}$$

Now we consider the fixed point $p_{n}$, the local coordinates of
the uniformizing system of $X_{({\bf g})}^{'}$ is
$[w_{i},\cdots,w_{n-1},1]$. Using the same method, we obtain the
$T$-equivariant Euler class of normal bundle of $p_{n} (i\leq
j\leq n)$ is given by
$$e_{T}(\nu_{p_{n}})=\prod_{k=1}^{n-i}\frac{d}{q_{n}}
\lambda_{k} \eqno{(5.17)}$$

Since $e_{1}^{*}(\eta_{1}),e_{2}^{*}(\eta_{2})$ and
$e_{3}^{*}(\eta_{3})$ all belong to $H^{*}(X_{({\bf g})}^{'},{\bf
Q})$. From the ordinary ring structure of weighted projective
space in section 5.1, we only consider $\xi_{1}\in H^{2}(X_{({\bf
g})}^{'},{\bf Q})$, the generator of $H^{2}(X_{({\bf g})}^{'},{\bf
Q})={\bf Q}$. Suppose $L\longrightarrow X_{({\bf g})}^{'}$ be the
canonical line bundle whose first chern-class is $\xi_{1}$. The
corresponding Cartier divisior is $D=D_{1}D_{2}\cdots
D_{i-1}D_{i}+\cdots+D_{1}D_{2}\cdots D_{i-1}D_{n}$, where
$D_{j}=\{z_{j}=0\}\subset {\bf P}(Q)$ is the basic divisor. Then
from Oda [20], in the neighbor $U_{j}\cap X_{({\bf g})}^{'}$,
$(i\leq j\leq n-1)$, let
$m=u_{1}\rho_{1}+\cdots+u_{n-i}\rho_{n-i}$, and
$<-m,v_{i}>=1$,$\cdots$, $<-m,v_{j-1}>=1$,
$<-m,v_{j+1}>=1$,$\cdots$,$<-m,v_{n}>=1$, then we calculate
$-m=\sum_{k\neq
r}\frac{d}{q_{n}}\rho_{k}-\frac{d}{q_{n}}(\sum_{k\neq
r}\frac{q_{k+i-1}}{q_{j}})\rho_{r}$, where $j=i+r-1$. So the
divisor $D$ is given by the rational function $\chi^{-m}$ on
$U_{j}\cap X_{({\bf g})}^{'}$. Similarly, it is given by the
rational function
$\chi^{-m}=\chi^{\frac{d}{q_{n}}(\rho_{1}+\cdots+\rho_{n-i})}$ on
$U_{n}\cap X_{({\bf g})}^{'}$. Hence the action of $T$ on the
corresponding line bundle of $D$ at the fixed points $p_{j}$ has
weights
$$\mbox{On}~~ p_{j},(i\leq j\leq n-1):  \sum_{k\neq
r}\frac{d}{q_{n}}\lambda_{k}-\frac{d}{q_{n}}\left(\sum_{k\neq
r}\frac{q_{k+i-1}}{q_{j}}\right)\lambda_{r} \eqno{(5.18)}$$
$$\mbox{On}~~ p_{n}:  \frac{d}{q_{n}}\sum_{k}\lambda_{k} \eqno{(5.19)}$$
On the other hand, we also can write $e_{1}^{*}\eta_{1}\wedge
e_{2}^{*}\eta_{2}\wedge e_{3}^{*}\eta_{3}= a(\xi_{1})^{s}$,
$a\in{\bf Q}$, $s$ is an integer.

Now we analyze the Euler form $e(E_{({\bf g})}^{''})$. From
section 5.5., we  compute the local generated vectors of the
obstruction bundle $E_{({\bf g})}$, and $E_{({\bf
g})}=\oplus_{l=1}^{e}E_{l}$. For each line bundle $E_{l}$, from
(5.11), we have the transition function of $E_{l}$ as
$$h_{jn}(x,c)=\left(x,x_{n}^{q_{t_{l}}/q_{n}}(x)\cdot c\right), ~~(i\leq j\leq n-1)$$
Because the line bundle $E_{l}^{''}$ is the reduction of $E_{l}$
under the ${\bf Z}_{d_{l}}$-invariant homomorphism, the transition
function of the line bundle $E_{l}^{''}$
$$(U_{n}\cap X_{({\bf g})}^{'})\times c\supset
(U_{n}\cap U_{j}\cap X_{({\bf g})}^{'})\times c\longrightarrow
(U_{n}\cap U_{j}\cap X_{({\bf g})}^{'})\times c\subset (U_{j}\cap
X_{({\bf g})}^{'})\times c$$ is given by
$$h_{jn}(x,c)=\left(x,(x_{n}^{q_{t_{l}}/q_{n}})^{d_{l}}(x)\cdot c\right),
~~(i\leq j\leq n-1)$$ So we can define the action of $T$ on
$E_{({\bf g})}^{''}=\bigoplus_{l=1}^{e}E_{l}^{''}$ as follows, on
the line bundle $E_{l}^{''}$:
$$
\begin{array}{ll} (1)&  t(x,c)=(tx,c)=\left(tx,\chi^{0}(t)c\right), t\in T,
(x,c)\in (U_{n}\cap X_{({\bf g})}^{'})\times {\bf C};\\
(2)& t(x,c)=\left(tx,(x_{n}^{q_{t_{l}}/q_{n}})^{d_{l}}(t)c\right)
=\left(tx,(\chi^{-\frac{d}{q_{j}}\rho_{r}})^{\frac{q_{t_{l}}d_{l}}{q_{n}}}(t)c\right),
\\
~~ & t\in T, (x,c)\in (U_{j}\cap X_{({\bf g})}^{'})\times {\bf C}.
\end{array}
$$
where $(i\leq j\leq n-1), j=i+r-1$. Then the action of $T$ on
$E_{l}^{''}$ at the fixed points $p_{n}, p_{j} (i\leq j\leq n-1)$
has weights
$$0,~~~ -\frac{q_{t_{l}}d_{l}}{q_{n}}\cdot \frac{d}{q_{j}}\lambda_{r}.
\eqno{(5.20)}$$

So from the localization formula, when we consider the orbifold
$X_{({\bf g})}^{'}$, see Corollary 9.13 in [12], we have the
integration $$
\int_{X_{(\bf{g})}^{'}}^{orb}e^{*}_{1}\eta_{1}\wedge
e^{*}_{2}\eta_{2}\wedge e^{*}_{3}\eta_{3}\wedge \Pi_{l=1}^{e}
e_{A_{l}^{''}}(E_{l}^{''})=
\frac{a\left(\frac{d}{q_{n}}\sum_{k=1}^{n-i}\lambda_{k}\right)^{s}\cdot
0 }{a_{n}\cdot\prod_{k=1}^{n-i}\frac{d}{q_{n}}\lambda_{k}} +$$
$$\sum_{j=i}^{n-1}\frac{a\left[\sum_{k\neq
r}\frac{d}{q_{n}}\lambda_{k}-\frac{d}{q_{n}}\left(\sum_{k\neq r
}\frac{q_{k+i-1}}{q_{j}}\right)\lambda_{r}\right]^{s}\cdot
\prod_{l=1}^{e}\left(-\frac{q_{t_{l}}d_{l}}{q_{n}}\cdot
\frac{d}{q_{j}}\lambda_{r}\right)}
{a_{j}\cdot\left(-\frac{d}{q_{j}}\lambda_{r}\right)\cdot
\prod_{k\neq
r}\frac{d}{q_{n}}\left(\lambda_{k}-\frac{q_{i+k-1}}{q_{j}}\lambda_{r}\right)}
\eqno{(5.21)}$$ Where $j=i+r-1$, $a_{j}$ is the order of the local
cyclic group of $p_{j}$ in the orbifold $X_{({\bf g})}^{'}$.


\subsection{Example.}

In this example we use the methods of the above sections to
calculate the 3-point functions. Let $Q=(1,2,2,3,3,3)$ and ${\bf
P}(Q)={\bf P}^{5}_{1,2,2,3,3,3}$ be the weighted projective space
of type $Q$, then $q_{0}=1, q_{1}=q_{2}=2,q_{3}=q_{4}=q_{5}=3$.
From proposition 3.2.2., we have $C_{0}=I_{5\times 5} $.  So let
$v_{1}=e_{1}, v_{2}=e_{2},
v_{3}=e_{3},v_{4}=e_{4},v_{5}=e_{5},v_{0}=-\Sigma\frac{q_{i}}{q_{0}}v_{i}=(-2,-2,-3,-3,-3)$.
The fan $\Xi$ of ${\bf P}_{1,2,2,3,3,3}^{5}$ is generated by
$\{v_{0},v_{1},v_{2},v_{3},v_{4},v_{5}\}$. For
$\sigma_{5}=(v_{0},v_{1},v_{2},v_{3},v_{4})$, we have
$G_{\sigma_{5}}=N/N_{\sigma_{5}}={\bf Z}_{3}$. We  write the
matrix representation of the action of ${\bf Z}_{3}$ on
$U_{\sigma_{5}^{'}}={\bf C}^{5}$ as follows.
$$
\left(
\begin{array}{ccccc}
1&0&0&0&0\\
0&1&0&0&0\\
0&0&1&0&0\\
0&0&0&1&0\\
0&0&0&0&1
\end{array}
\right); \left (
\begin{array}{ccccc}
e^{2\pi i\cdot \frac{1}{3}}&0&0&0&0\\
0&e^{2\pi i\cdot \frac{2}{3}}&0&0&0\\
0&0&e^{2\pi i\cdot \frac{2}{3}}&0&0\\
0&0&0&1&0\\
0&0&0&0&1
\end{array}
\right)
$$
$$
\left (
\begin{array}{ccccc}
e^{2\pi i\cdot \frac{2}{3}}&0&0&0&0\\
0&e^{2\pi i\cdot \frac{1}{3}}&0&0&0\\
0&0&e^{2\pi i\cdot \frac{1}{3}}&0&0\\
0&0&0&1&0\\
0&0&0&0&1
\end{array}
\right)
$$
For $\sigma_{4}=(v_{0},v_{1},v_{2},v_{3},v_{5})$, and
$\sigma_{3}=(v_{0},v_{1},v_{2},v_{4},v_{5})$, we have
$G_{\sigma_{4}}=G_{\sigma_{3}}={\bf Z}_{3}$.  The actions of ${\bf
Z}_{3}$ on $U_{\sigma_{4}^{'}}=U_{\sigma_{3}^{'}}={\bf C}^{5}$ are
the same  as above.

For $\sigma_{2}=(v_{0},v_{1},v_{3},v_{4},v_{5})$. We  write the
matrix representation of the action of ${\bf Z}_{2}$ on
$U_{\sigma_{2}^{'}}={\bf C}^{5}$ as follows.
$$
\left(
\begin{array}{ccccc}
1&0&0&0&0\\
0&1&0&0&0\\
0&0&1&0&0\\
0&0&0&1&0\\
0&0&0&0&1
\end{array}
\right); \left (
\begin{array}{ccccc}
e^{2\pi i\cdot \frac{1}{2}}&0&0&0&0\\
0&1&0&0&0\\
0&0&e^{2\pi i\cdot \frac{1}{2}}&0&0\\
0&0&0&e^{2\pi i\cdot \frac{1}{2}}&0\\
0&0&0&0&e^{2\pi i\cdot \frac{1}{2}}
\end{array}
\right)
$$
For $\sigma_{1}=(v_{0},v_{2},v_{3},v_{4},v_{5})$, we have
$G_{\sigma_{1}}={\bf Z}_{2}$.  The actions of ${\bf Z}_{2}$ on
$U_{\sigma_{1}^{'}}={\bf C}^{5}$ is the same as above.

For $\sigma_{0}=(v_{1},v_{2},v_{3},v_{4},v_{5})$,
$G_{\sigma_{0}}=1$, and the action is trivial. If we let
$$
g_{1}=\left (
\begin{array}{ccccc}
e^{2\pi i\cdot \frac{1}{3}}&0&0&0&0\\
0&e^{2\pi i\cdot \frac{2}{3}}&0&0&0\\
0&0&e^{2\pi i\cdot \frac{2}{3}}&0&0\\
0&0&0&1&0\\
0&0&0&0&1
\end{array}
\right); g_{2}= \left (
\begin{array}{ccccc}
e^{2\pi i\cdot \frac{1}{2}}&0&0&0&0\\
0&1&0&0&0\\
0&0&e^{2\pi i\cdot \frac{1}{2}}&0&0\\
0&0&0&e^{2\pi i\cdot \frac{1}{2}}&0\\
0&0&0&0&e^{2\pi i\cdot \frac{1}{2}}
\end{array}
\right)
$$
Then we have the twisted  sectors:
$X_{(g_{1})}=X_{(g_{1}^{2})}={\bf P}(Q_{\tau})$, where
$Q_{\tau}=(0,0,0,3,3,3)$, $\tau=(v_{0},v_{1},v_{2})$;
$X_{(g_{2})}={\bf P}(Q_{\delta})$, $Q_{\delta}=(0,2,2,0,0,0),
\delta=(v_{0},v_{3},v_{4},v_{5})$. The degree shifting numbers:
$\iota_{(g_{1})}=\frac{5}{3},\iota_{(g_{1}^{2})}=\frac{4}{3},
\iota_{(g_{2})}=2$. So the Chen-Ruan  cohomology group of ${\bf
P}^{5}_{1,2,2,3,3,3}$ is:
\begin{eqnarray}
H^{d}_{orb}({\bf P}^{5}_{1,2,2,3,3,3};{\bf Q})&=&H^{d}({\bf P}^{5}_{1,2,2,3,3,3};{\bf Q})\nonumber \\
&\oplus & H^{d-\frac{10}{3}}({\bf P}(Q_{\tau});{\bf Q}) \nonumber \\
&\oplus & H^{d-\frac{8}{3}}({\bf P}(Q_{\tau});{\bf Q}) \nonumber \\
&\oplus & H^{d-4}({\bf P}(Q_{\delta});{\bf Q}) \nonumber
\end{eqnarray}

All the 3-multisectors are:
$X_{(g_{1},g_{1},g_{1})}=X_{(g_{1}^{2},g_{1}^{2},g_{1}^{2})}={\bf
P}(Q_{\tau})$, $X_{(g_{1},g_{1}^{2},1)}={\bf P}(Q_{\tau})$,
$X_{(g_{2},g_{2},1)}={\bf P}(Q_{\delta})$. In the 3-multisectors
$X_{(g_{1},g_{1}^{2},1)}$ and $X_{(g_{2},g_{2},1)}$, from (5.2),
the dimension of the  obstruction bundle of these two
3-multisectors are all zero, so the integration (5.3) is the usual
integration on orbifold. The orbifold cup product can be described
easily.

For the 3-multisector $X_{(g_{1}^{2},g_{1}^{2},g_{1}^{2})}$, the
dimension of the obstruction bundle $E_{({\bf g})}$ is 1. Let
$X_{({\bf g})}=X_{(g_{1}^{2},g_{1}^{2},g_{1}^{2})}$, $\eta_{j}\in
H^{*}(X_{(g_{1}^{2})};{\bf Q}), (j=1,2,3)$, then
$$<\eta_{1},\eta_{2},\eta_{3}>_{orb}=\int_{X_{(\bf{g})}}^{orb}e^{*}_{1}\eta_{1}\wedge
e^{*}_{2}\eta_{2}\wedge e^{*}_{3}\eta_{3}\wedge
e_{A}(E_{(\bf{g})}) \eqno{(5.22)}$$

Next we use the localization formula (5.21) to compute the 3-point
function (5.22). First we describe the obstruction bundle
$E_{({\bf g})}$ over $X_{({\bf g})}$. We see that $g_{1}^{2}$
generates the cyclic group ${\bf Z}_{3}$, so $K({\bf g})={\bf
Z}_{3}$. To describe the obstruction bundle, we consider the
orbifold sphere $(S^{2},(x_{1},x_{2},x_{3}),(3,3,3))$,  the
orbifold fundamental group is:
$$\pi_{1}^{orb}(S^{2})=\left\{\lambda_{1},\lambda_{2},\lambda_{3}|
\lambda_{i}^{3}=1, \lambda_{1}\lambda_{2}\lambda_{3}=1\right\}.$$
Its orbifold universal cover is the Euclidean plane $E^{2}$([24]).
We use the triangle group model $\triangle^{*}(p,q,r)$ for
$\frac{1}{p}+\frac{1}{q}+\frac{1}{r}=1$. Then from [24], the full
triangle group $\triangle^{*}(p,q,r)$ is defined to be the
triangle group of isometries of $E^{2}$ generated by the
reflexions $L$,$M$ and $N$ in the three sides $YZ$,$ZX$ and $XY$
of $\triangle$. In each case, it is easy to see that the
translates of $\triangle$ by $\triangle^{*}(3,3,3)$ title $E^{2}$.
This tiling is shown in [24].

Now $\triangle^{*}(3,3,3)$ has a natural subgroup of index two-the
orientation preserving subgroup which is denoted
$\triangle(3,3,3)$. The product $LM$ of two of the generating
reflexions of $\triangle^{*}(3,3,3)$ is a rotation through
$2\pi/3$ about $Z$ and $MN$ and $NL$ are also rotations. Write
$LM=\lambda_{1}$, $MN=\lambda_{2}$, $NL=\lambda_{3}$ so that
$\lambda_{1},\lambda_{2},\lambda_{3}$ are rotations about $Z,X,Y$
respectively, see Figure.1.  Clearly
$\lambda_{1},\lambda_{2},\lambda_{3}$ lie in $\triangle(3,3,3)$
and one can show that they generate it. So
$\triangle(3,3,3)=\pi_{1}^{orb}(S^{2})$, and
$E^{2}/\pi_{1}^{orb}(S^{2})=S^{2}(3,3,3)$. Let $P=\triangle\cup
L\triangle$, then the $E^{2}/\pi_{1}^{orb}(S^{2})$ is obtained
from the fundamental region $\in=P$ by the following manner: $XZ$
identified with $LXZ$, and $XY$ identified with $LXY$.

Consider the homomorphism $\rho:
\pi_{1}^{orb}(S^{2})\longrightarrow K({\bf g})={\bf Z}_{3}$ given
by $\lambda_{i}\longmapsto g_{1}^{2}$. $ker(\rho)$ is a normal
subgroup generated by commutators of $\lambda_{i}$ and the element
$\lambda_{1}\lambda_{2}^{-1}$. $ker(\rho)$ acts freely on $E^{2}$
with quotient being the hexigon shown in Figure.2.

\begin{center}
\includegraphics[scale=.7]{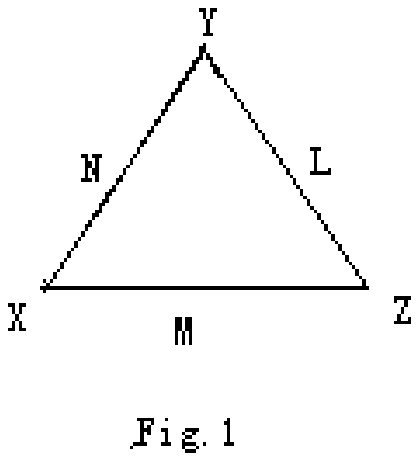}
\includegraphics[scale=.7]{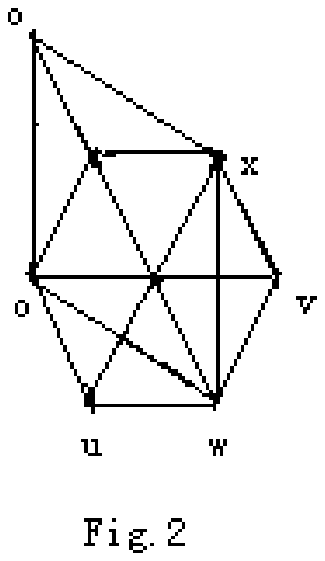}
\end{center}

Suppose $\alpha=\overrightarrow{ouw}$,
$\beta=\overrightarrow{wvx}$. Then we obtain a symplectic basis
$\{a,b\}=\{\alpha,\beta\}$ for $H_{1}(\Sigma;{\bf Z})$. Taking the
Universal Coefficients Theorem  dual, we obtain a canonical
symplectic basis $\{\check{\alpha},\check{\beta}\}$ for
$H^{1}(\Sigma;{\bf Z})$. From Riemann bilinear relation in [15],
we know that there exists a complex basis $\omega$ for
$H^{1,0}(\Sigma)$ such that its period matrix with respect to the
$a$-class is the identity matrix and its period matrix with
respect to the $b$-class is a symmetric complex matrix $R=(r)$
with $ImR$ a positive matrix. Here we can write
$\omega=\check{\alpha}+r\check{\beta}$.

The rotation $\lambda_{1}$ induces the following map on
$\pi_{1}(\Sigma)$:
$$
\left\{
\begin{array}{l}
\alpha\longmapsto\beta^{-1};\\
\beta\longmapsto\alpha\beta^{-1}.
\end{array}
\right.
$$
Hence the automorphism $(\lambda_{1})_{*}:
H_{1}(\Sigma)\longrightarrow H_{1}(\Sigma)$ can be expressed in
the matrix:
$$
\left(
\begin{array}{cc}
0&1\\
-1&-1
\end{array}
\right)
$$
with respect to the basis $\{\alpha,\beta\}$. The automorphism
$\lambda_{1}^{*}: H^{1}(\Sigma)\longrightarrow H^{1}(\Sigma)$ can
be expressed by the matrix:
$$
A=\left(
\begin{array}{cc}
0&-1\\
1&-1
\end{array}
\right)
$$
with respect to the dual basis $\{\check{\alpha},\check{\beta}\}$,
which is the transpose of the above matrix.

From the matrix $A$, we  calculate that
$\lambda_{1}^{*}\omega=-r\check{\alpha}+(-r+1)\check{\beta}$.
Since $\lambda_{1}$ is a holomorphic map, $\lambda_{1}^{*}$
preserves the subspace $H^{1,0}(\Sigma)$. Hence
$\lambda_{1}^{*}\omega$ can be expressed as complex linear
combination of $\omega$, so we have
$\lambda_{1}^{*}\omega=-r\omega$. We obtain $-r^{2}=-r+1$,
$r=\frac{1\pm \sqrt{3}i}{2}$. While $ImR>0$, $r=\frac{1+
\sqrt{3}i}{2}$. Because
$H^{0,1}(\Sigma)=\overline{H^{1,0}(\Sigma)}$, $H^{0,1}(\Sigma)$ is
generated by $\overline{\omega}$,
$\overline{\omega}=\check{\alpha}+\bar{r}\check{\beta}$

Now let $\tau=(v_{0},v_{1},v_{2})$, then
$\overline{O}_{\tau}=X_{({\bf g})}$. Suppose $U_{i}=\{z_{i}\neq
0\}$, then $X_{({\bf g})}=(X_{({\bf g})}\cap U_{3})\cup (X_{({\bf
g})}\cap U_{4})\cup (X_{({\bf g})}\cap U_{5})$. Let
$U_{3}=V_{3}/{\bf Z}_{3}$, $U_{4}=V_{4}/{\bf Z}_{3}$,
$U_{5}=V_{5}/{\bf Z}_{3}$, from section 5.5., the coordinates of
$V_{3}$,$V_{4}$ and $V_{5}$ are:
$$
\begin{array}{ll}
V_{5}: &
\left\{x_{0}=\frac{z_{0}}{(z_{5})^{1/3}},x_{1}=\frac{z_{1}}{(z_{5})^{2/3}},x_{2}=\frac{z_{2}}{(z_{5})^{2/3}},
x_{3}=\frac{z_{3}}{z_{5}},x_{4}=\frac{z_{4}}{z_{5}},x_{5}=1\right\};\\
V_{4}:&
\left\{y_{0}=\frac{z_{0}}{(z_{4})^{1/3}},y_{1}=\frac{z_{1}}{(z_{4})^{2/3}},y_{2}=\frac{z_{2}}{(z_{4})^{2/3}},
y_{3}=\frac{z_{3}}{z_{4}},y_{4}=1,y_{5}=\frac{z_{5}}{z_{4}}\right\};\\
V_{3}:&
\left\{w_{0}=\frac{z_{0}}{(z_{3})^{1/3}},w_{1}=\frac{z_{1}}{(z_{3})^{2/3}},w_{2}=\frac{z_{2}}{(z_{3})^{2/3}},
w_{3}=1,w_{4}=\frac{z_{4}}{z_{3}},w_{5}=\frac{z_{5}}{z_{3}}\right\}.
\end{array}
$$
For the chart $(V_{5},{\bf Z}_{3},\pi_{5})$, $TV_{5}|_{p_{5}}$ has
framing  $\{\frac{\partial}{\partial
x_{0}},\frac{\partial}{\partial x_{1}},\frac{\partial}{\partial
x_{2}},\frac{\partial}{\partial x_{3}},\frac{\partial}{\partial
x_{4}}\}$. So $TV_{5}|_{p_{5}}\otimes H^{0,1}(\Sigma)$  has
framing $$\left\{\frac{\partial}{\partial x_{0}}\otimes
\overline{\omega},\frac{\partial}{\partial x_{1}}\otimes
\overline{\omega},\frac{\partial}{\partial x_{2}}\otimes
\overline{\omega},\frac{\partial}{\partial x_{3}}\otimes
\overline{\omega},\frac{\partial}{\partial x_{4}}\otimes
\overline{\omega}\right\}$$ Because
$\lambda_{1}^{*}\overline{\omega}=-\bar{r}\overline{\omega}$, and
the action of $\lambda_{1}^{*}$ on $TV_{5}|_{p_{5}}$ is expressed
by the matrix $g_{1}^{2}$, so the representation matrix of the
action of  $\lambda_{1}^{*}$ on $TV_{5}|_{p_{5}}\otimes
H^{0,1}(\Sigma)$ is:
$$
(\ast)\left(
\begin{array}{ccccc}
1&0&0&0&0\\
0&-\frac{1}{2}-\frac{\sqrt{3}}{2}i&0&0&0\\
0&0&-\frac{1}{2}-\frac{\sqrt{3}}{2}i&0&0\\
0&0&0&-\frac{1}{2}+\frac{\sqrt{3}}{2}i&0\\
0&0&0&0&-\frac{1}{2}+\frac{\sqrt{3}}{2}i
\end{array}
\right)
$$
So the matrix $(\ast)$ has eigenvalue 1 with multiplication one,
and the corresponding eigenvector is $(1,0,0,0,0)$, the generator
of $(TV_{5}|_{p_{5}}\otimes H^{0,1}(\Sigma))^{K({\bf g})}$ is
$\xi_{0}\otimes \overline{\omega}=\frac{\partial}{\partial
x_{0}}\otimes \overline{\omega}$. Similarly, we  compute the
generator of $(TV_{4}|_{p_{4}}\otimes H^{0,1}(\Sigma))^{K({\bf
g})}$ is $\xi_{0}^{'}\otimes
\overline{\omega}=\frac{\partial}{\partial y_{0}}\otimes
\overline{\omega}$ and the generator of $(TV_{3}|_{p_{3}}\otimes
H^{0,1}(\Sigma))^{K({\bf g})}$ is $\xi_{0}^{''}\otimes
\overline{\omega} =\frac{\partial}{\partial w_{0}}\otimes
\overline{\omega}$.

Now we describe the local uniformizing charts for $E_{({\bf g})}$:

$\bullet$ $\forall x\in X_{({\bf g})}$, then $C({\bf
g})=G_{x}=K({\bf g})={\bf Z}_{3}$, and $(V_{x}^{{\bf g}}\times
{\bf C},K({\bf g}),\widetilde{\pi})$ is a uniformizing system for
$E_{({\bf g})}$ where $K({\bf g})$ acts on $V_{x}^{{\bf g}}\times
{\bf C}$ by: $g_{1}^{2}(u,v)=(u,e^{2\pi i\cdot \frac{2}{3}}v)$.

The bundle $E_{({\bf g})}$ is a line bundle, $dim_{{\bf
C}}X_{({\bf g})}=2$, so the 3-point function (5.22) is nonzero
only if there is some $\eta_{i}\in H^{2}(X_{(g_{1})};{\bf Q})$.
Without loss of generality, assume $\eta_{1}\in
H^{2}(X_{(g_{1}^{2})};{\bf Q})$,  $\eta_{2}\in
H^{0}(X_{(g_{1}^{2})};{\bf Q})$, $\eta_{3}\in
H^{0}(X_{(g_{1}^{2})};{\bf Q})$. In this case
$$<\eta_{1},\eta_{2},\eta_{3}>_{orb}=\eta_{2}\eta_{3}\int_{X_{(\bf{g})}}^{orb}\eta_{1}
\wedge e_{A}(E_{(\bf{g})}) \eqno{(5.23)}$$

From the first part of this section, we see that in this case, the
orbifold principal $S^{1}$ bundle is $P_{({\bf g})}$, and let
$E_{({\bf g})}^{'}=(P_{({\bf g})}/K({\bf g}))\times_{S^{1}}{\bf
C}$ over $X_{({\bf g})}$. $\pi_{K({\bf g})}:P_{({\bf
g})}\longrightarrow P_{({\bf g})}/K({\bf g})$. Note that
$\pi_{K({\bf g})}$ on each fibre is given by $z\longmapsto z^{3}$.
The lie algebra of $F=S^{1}$ can be identified with ${\bf R}$.
Hence the induced map on the lie algebra $(\pi_{K({\bf g})})_{*}:
{\bf R}\longrightarrow {\bf R} $ is just multiplication by 3, so
from (5.14)
$$\int_{X_{(\bf{g})}}^{orb}\eta_{1}
\wedge
e_{A}(E_{(\bf{g})})=\frac{1}{3}\int_{X_{(\bf{g})}}^{orb}\eta_{1}
\wedge e_{A^{'}}(E_{(\bf{g})}^{'}) $$ Where $A$ and $A^{'}$ is the
connections of $E_{(\bf{g})}$ and $E_{(\bf{g})}^{'}$ such that
$\pi_{K({\bf g})}^{*}(A^{'})=A$. Then $E_{({\bf g})}^{'}$ induces
an orbifold bundle $E_{(\bf{g})}^{''}$ over the reduced orbifold
$X_{({\bf g})}^{'}$, from (5.15)
$$\int_{X_{(\bf{g})}}^{orb}\eta_{1}
\wedge
e_{A^{'}}(E_{(\bf{g})}^{'})=\frac{1}{3}\int_{X_{(\bf{g})}^{'}}^{orb}\eta_{1}
\wedge e_{A^{''}}(E_{(\bf{g})}^{''}) $$ Where $A^{''}$ is the
connection of $E_{(\bf{g})}^{''}$ induced  from
$E_{(\bf{g})}^{'}$. So we obtain
$$\int_{X_{(\bf{g})}}^{orb}\eta_{1}
\wedge
e_{A}(E_{(\bf{g})})=\frac{1}{9}\int_{X_{(\bf{g})}^{'}}^{orb}\eta_{1}
\wedge e_{A^{''}}(E_{(\bf{g})}^{''}) \eqno{(5.24)}$$

We now compute the integration
$\int_{X_{(\bf{g})}^{'}}^{orb}\eta_{1} \wedge
e_{A^{''}}(E_{(\bf{g})}^{''})$ in (5.24). The uniformizing system
of $E_{({\bf g})}^{''}$ over $X_{({\bf g})}^{'}$ can be described
as follows.

$\bullet$ $\forall x\in X_{({\bf g})}^{'}$, then $C({\bf
g})/K({\bf g})=1$ is the trivial group.

Now we use the localization technique  to calculate the
integration (5.24). Note that $X_{({\bf g})}^{'}=O_{\tau}$,
$\tau=(v_{0},v_{1},v_{2})$, then $N(\tau)=N/N_{\tau}$, where
$N_{\tau}$ is the sublattice generated by $\tau$, and
$M(\tau)=\tau^{\perp}\cap M$. The 2-torus associated to $X_{({\bf
g})}^{'}$ is $T=spec({\bf C}[M(\tau)])=O_{\tau}$. The characters
$\chi^{m}$ correspond to rational functions on $X_{({\bf g})}^{'}$
when $m\in M(\tau)$. Let $\{m_{1},m_{2},m_{3},m_{4},m_{5}\}$ be
the standard basis of $M$. Then $\{\rho_{1}=m_{3}-m_{5},
\rho_{2}=m_{4}-m_{5}\}$ is a basis for $M(\tau)$. The $T$-action
on $X_{({\bf g})}^{'}$ has three fixed points $p_{3},p_{4},p_{5}$.
First we study the action of $T$ on the normal bundle of
$p_{3},p_{4}$ and $p_{5}$, i.e., the orbifold tangent space of
$p_{3},p_{4}$ and $p_{5}$.

From (5.16) and (5.17), of course we can compute using the same
method as in section 5.6, we see that the $T$-equivariant Euler
class of the normal bundle of $p_{5}$ is given by
$e_{T}\left(\nu_{p_{5}}\right)=\lambda_{1}\lambda_{2}$, and we
also have:
$e_{T}\left(\nu_{p_{4}}\right)=(\lambda_{1}-\lambda_{2})(-\lambda_{2})$,
$e_{T}\left(\nu_{p_{3}}\right)=(-\lambda_{1}+\lambda_{2})(-\lambda_{1})$.
In particular, we have $\chi^{-\rho_{2}}=y_{5}$ in the
neighborhood $V_{4}$, and $\chi^{-\rho_{1}}=\omega_{5}$ in the
neighborhood $V_{5}$.

The orbifold line bundle $E_{({\bf g})}^{''}$ is trivialized by
the generator $\frac{\partial}{\partial x_{0}}\otimes
\overline{\omega}$ on $U_{5}\cap X_{({\bf g})}^{'}$,
$\frac{\partial}{\partial y_{0}}\otimes \overline{\omega}$ on
$U_{4}\cap X_{({\bf g})}^{'}$, and $\frac{\partial}{\partial
w_{0}}\otimes \overline{\omega}$ on $U_{3}\cap X_{({\bf g})}^{'}$.
On $V_{4}$: $\frac{\partial}{\partial x_{0}}=\sum \frac{\partial
y_{j}}{\partial x_{0}}\frac{\partial}{\partial
y_{j}}=\frac{\partial}{\partial
x_{0}}(\frac{x_{0}}{x_{4}^{1/3}})\frac{\partial}{\partial
y_{0}}=y_{5}^{1/3}\frac{\partial}{\partial y_{0}}$; On $V_{3}$:
$\frac{\partial}{\partial x_{0}}=\sum \frac{\partial
w_{j}}{\partial x_{0}}\frac{\partial}{\partial
w_{j}}=\frac{\partial}{\partial
x_{0}}(\frac{x_{0}}{x_{3}^{1/3}})\frac{\partial}{\partial
w_{1}}=w_{5}^{1/3}\frac{\partial}{\partial w_{0}}$. So the
transition function:
$$h_{45}:(U_{5}\cap X_{({\bf g})}^{'})\times {\bf C}\supset
 (U_{4}\cap U_{5}\cap X_{({\bf g})}^{'})\times {\bf C}\longrightarrow
 (U_{4}\cap U_{5}\cap X_{({\bf g})}^{'})\times {\bf C}\subset
(U_{4}\cap X_{({\bf g})}^{'})\times {\bf C}$$ can be written as:
 $$h_{45}(x,c)=\left(x,(y_{5}^{1/3})^{3}(x)\cdot c\right).$$
The transition function
$$h_{35}:(U_{5}\cap X_{({\bf
g})}^{'})\times {\bf C}\supset
 (U_{5}\cap U_{3}\cap X_{({\bf g})}^{'})\times {\bf C}\longrightarrow
 (U_{5}\cap U_{3}\cap X_{({\bf g})}^{'})\times {\bf C}\subset
(U_{3}\cap X_{({\bf g})}^{'})\times {\bf C}$$ is
 $$h_{35}(x,c)=\left(x,(w_{5}^{1/3})^{3}(x)\cdot c\right)$$
So we can define the $T$ action on $E_{({\bf g})}^{''}$ by
$$
\begin{array}{ll} (1)&  t(x,c)=(tx,c)=(tx,\chi^{0}(t)c), t\in T,
(x,c)\in (U_{5}\cap X_{({\bf g})}^{'})\times {\bf C};\\
(2)&  t(x,c)=(tx,y_{5}(t)c)=(tx,\chi^{-\rho_{2}}(t)c), t\in T,
(x,c)\in (U_{4}\cap X_{({\bf g})}^{'})\times {\bf C};\\
(3)&  t(x,c)=(tx,w_{5}(t)c)=(tx,\chi^{-\rho_{1}}(t)c), t\in T,
(x,c)\in (U_{3}\cap X_{({\bf g})}^{'})\times {\bf C}.
\end{array}
$$
Then the  action of $T$ on $E_{({\bf g})}^{''}$ at the fixed
points $p_{5},p_{4},p_{3}$ has weights $0,(-\lambda_{2}),
(-\lambda_{1})$ respectively. we see that this can be calculated
from (5.18) and (5.19).

Since $\eta_{1}\in H^{2}(X_{(g_{1}^{2})};{\bf Q})$, so we can
associate $\eta_{1}$ with a weil divisor
$\{U_{\sigma},\chi^{-m^{\sigma}}\}$. As before the transition
functions $h_{\tau\sigma}: U_{\sigma}\times{\bf C}\supset
U_{\sigma\cap\tau}\times {\bf C}\longrightarrow
U_{\sigma\cap\tau}\times {\bf C}\subset U_{\tau}\times {\bf C}$
for the line bundle are given by
$h_{\tau\sigma}(x,c)=(x,\chi^{(m_{\sigma}-m_{\tau})}(x)c)$. Then
one can define a $T$- action on the bundle that makes it a
$T$-equivariant bundle ([20]) as
follows:$t(x,c)=(tx,\chi^{-m_{\sigma}}(t)c)$ for $t\in T$ and
$(x,c)\in U_{\sigma}\times {\bf C}$.

Suppose $D_{i}=\{z_{i}=0\}(i=0,1,2,3,4,5)$ are the divisors of
$X$. We  see that the divisor of $X_{({\bf g})}^{'}$ can be
expressed as the combination of
$D_{0}D_{1}D_{2}D_{3},D_{0}D_{1}D_{2}D_{4},
\\ D_{0}D_{1}D_{2}D_{5}$. We consider the general $\eta_{1}$, and
suppose it can be expressed as
$D=a_{3}D_{0}D_{1}D_{2}D_{3}+a_{4}D_{0}D_{1}D_{2}D_{4}+a_{5}D_{0}D_{1}D_{2}D_{5}$.
Then in $U_{5}\cap X_{({\bf g})}^{'}$, suppose
$m=u_{1}\rho_{1}+u_{2}\rho_{2}$, let $<-m,v_{3}>=a_{3},
<-m,v_{4}>=a_{4}$, then $u_{1}=-a_{3}, u_{2}=-a_{4}$,
$\chi^{-m}=\chi^{a_{3}\rho_{1}+a_{4}\rho_{2}}$. The divisor $D$ is
given by the rational function
$\chi^{a_{3}\rho_{1}+a_{4}\rho_{2}}$ on $U_{5}\cap X_{({\bf
g})}^{'}$. Similarly, it is given by the rational function
$\chi^{-m}=\chi^{a_{3}\rho_{1}-(a_{5}+a_{3})\rho_{2}}$ on
$U_{4}\cap X_{({\bf g})}^{'}$ and
$\chi^{-(a_{5}+a_{4})\rho_{1}+a_{4}\rho_{2}}$ on $U_{3}\cap
X_{({\bf g})}^{'}$. Hence the action of $T$ on the corresponding
line bundle of $\eta_{1}$ at the fixed points $p_{5},p_{4}$ and
$p_{3}$ has weights $a_{3}\lambda_{1}+a_{4}\lambda_{2},
a_{3}\lambda_{1}-(a_{5}+a_{3})\lambda_{2}$ and
$-(a_{5}+a_{4})\lambda_{1}+a_{4}\lambda_{2}$ respectively. So
using the localization formula (5.21), we have:
\begin{eqnarray}
\int_{X_{(\bf{g})}^{'}}^{orb}\eta_{1} \wedge
e_{A^{''}}(E_{(\bf{g})}^{''})&=&\frac{(a_{3}\lambda_{1}+a_{4}\lambda_{2})\cdot
0}{\lambda_{1}\lambda_{2}}+\frac{(a_{3}\lambda_{1}-(a_{5}+a_{3})\lambda_{2})(-\lambda_{2})}
{(\lambda_{1}-\lambda_{2})(-\lambda_{2})} \nonumber \\ &+&
\frac{(-(a_{5}+a_{4})\lambda_{1}+a_{4}\lambda_{2})(-\lambda_{1})}
{(-\lambda_{1}+\lambda_{2})(-\lambda_{1})} \nonumber \\
&=& (a_{3}+a_{4}+a_{5}) \nonumber
\end{eqnarray}
So from (5.24), we have:
 $$\int_{X_{({\bf
g})}}^{orb}\eta_{1}\wedge e_{A}(E_{({\bf
g})})=\frac{1}{9}(a_{3}+a_{4}+a_{5}).$$ From (5.23),
$$<\eta_{1},\eta_{2},\eta_{3}>_{orb}=\frac{1}{9}(a_{3}+a_{4}+a_{5})\eta_{2}\eta_{3}.$$
For example, if $D=D_{0}D_{1}D_{2}D_{3}$, then $a_{3}=1,
a_{4}=a_{5}=0$, so
$$<\eta_{1},\eta_{2},\eta_{3}>_{orb}=\frac{1}{9}\eta_{2}\eta_{3}.$$

For the 3-multisector $X_{(g_{1},g_{1},g_{1})}={\bf P}(Q_{\tau})$,
the dimension of the obstruction bundle $E_{({\bf g})}$ is 2. Let
$X_{({\bf g})}=X_{(g_{1},g_{1},g_{1})}$, $\eta_{j}\in
H^{*}(X_{(g_{1})};{\bf Q}), (j=1,2,3)$, then
$$<\eta_{1},\eta_{2},\eta_{3}>_{orb}=\int_{X_{(\bf{g})}}^{orb}e^{*}_{1}\eta_{1}\wedge
e^{*}_{2}\eta_{2}\wedge e^{*}_{3}\eta_{3}\wedge
e_{A}(E_{(\bf{g})}) \eqno{(5.25)}$$

Next we use the localization technique to calculate the 3-point
function (5.25). First we describe the obstruction bundle
$E_{({\bf g})}$ over $X_{({\bf g})}$. We know that $g_{1}$
generate the group ${\bf Z}_{3}$, so $K({\bf g})={\bf Z}_{3}$. We
consider the orbifold sphere
$(S^{2},(x_{1},x_{2},x_{3}),(3,3,3))$, let $$\rho:
\pi_{1}^{orb}(S^{2})\longrightarrow K({\bf g})={\bf Z}_{3}$$is the
surjective homomorphism defined by $\lambda_{i}\longmapsto g_{1}$.
Then we have a Riemann surface $\Sigma=E^{2}/ker(\rho)$. Using the
same method  above, we have a basis  $\overline{\omega}$ of
$H^{0,1}(\Sigma)$, and
$$\lambda_{1}^{*}(\overline{\omega})=-\bar{r}\overline{\omega}$$

Let $\tau=(v_{0},v_{1},v_{2})$, then we know that
$\overline{O}_{\tau}=X_{({\bf g})}$. Assume $U_{i}=\{z_{i}\neq
0\}$, then $X_{({\bf g})}=(X_{({\bf g})}\cap U_{3})\cup (X_{({\bf
g})}\cap U_{4})\cup (X_{({\bf g})}\cap U_{5})$. Write
$U_{3}=V_{3}/{\bf Z}_{3}$, $U_{4}=V_{4}/{\bf Z}_{3}$,
$U_{5}=V_{5}/{\bf Z}_{3}$, then we can choose the coordinates of
the open set $V_{3}$, $V_{4}$ and $V_{5}$ the same as before.

For the coordinate neighborhood $(V_{5},{\bf Z}_{3},\pi_{5})$,
$TV_{5}|_{p_{5}}$ has a framing:  $\{\frac{\partial}{\partial
x_{0}},\frac{\partial}{\partial x_{1}},\frac{\partial}{\partial
x_{2}}, \\ \frac{\partial}{\partial x_{3}},
\frac{\partial}{\partial x_{4}}\}$, so $TV_{5}|_{p_{5}}\otimes
H^{0,1}(\Sigma)$ has a framing
$$\left\{\frac{\partial}{\partial x_{0}}\otimes
\overline{\omega},\frac{\partial}{\partial x_{1}}\otimes
\overline{\omega},\frac{\partial}{\partial x_{2}}\otimes
\overline{\omega},\frac{\partial}{\partial x_{3}}\otimes
\overline{\omega},\frac{\partial}{\partial x_{4}}\otimes
\overline{\omega}\right\}$$ Because
$\lambda_{1}^{*}\overline{\omega}=-\bar{r}\overline{\omega}$, and
the action $\lambda_{1}^{*}$ on $TV_{5}|_{p_{5}}$ can be
represented by the matrix $g_{1}$, so the representation matrix of
the action of $\lambda_{1}^{*}$ on $TV_{5}|_{p_{5}}\otimes
H^{0,1}(\Sigma)$ is:
$$
(\ast \ast)\left(
\begin{array}{ccccc}
-\frac{1}{2}-\frac{\sqrt{3}}{2}i&0&0&0&0\\
0&1&0&0&0\\
0&0&1&0&0\\
0&0&0&-\frac{1}{2}+\frac{\sqrt{3}}{2}i&0\\
0&0&0&0&-\frac{1}{2}+\frac{\sqrt{3}}{2}i
\end{array}
\right)
$$

We see that the matrix $(\ast \ast)$ has eigenvalue 1 with
multiplication 2, the corresponding eigenvectors are:
$(0,1,0,0,0)$ and $(0,0,1,0,0)$, so the invariant subspace
$(TV_{5}|_{p_{5}}\otimes H^{0,1}(\Sigma))^{K({\bf g})}$ has
generators: $\xi_{1}\otimes
\overline{\omega}=\frac{\partial}{\partial x_{1}}\otimes
\overline{\omega}$ and $\xi_{2}\otimes
\overline{\omega}=\frac{\partial}{\partial x_{2}}\otimes
\overline{\omega}$. Similarly, $(TV_{4}|_{p_{4}}\otimes
H^{0,1}(\Sigma))^{K({\bf g})}$ has generators: $\xi_{1}^{'}\otimes
\overline{\omega}=\frac{\partial}{\partial y_{1}}\otimes
\overline{\omega}$ and $\xi_{2}^{'}\otimes
\overline{\omega}=\frac{\partial}{\partial y_{2}}\otimes
\overline{\omega}$; $(TV_{3}|_{p_{3}}\otimes
H^{0,1}(\Sigma))^{K({\bf g})}$ has generators:
$\xi_{1}^{''}\otimes \overline{\omega} =\frac{\partial}{\partial
w_{1}}\otimes \overline{\omega}$ and $\xi_{2}^{''}\otimes
\overline{\omega} =\frac{\partial}{\partial w_{2}}\otimes
\overline{\omega}$.

We describe the uniformizing system of $E_{({\bf g})}$ as follows.

$\bullet$ $\forall x\in X_{({\bf g})}$, then $C({\bf
g})=G_{x}=K({\bf g})={\bf Z}_{3}$, and if $(V_{x}^{{\bf g}}\times
{\bf C}^{2},K({\bf g}),\widetilde{\pi})$ is a uniformizing system
of the bundle $E_{({\bf g})}$, then  $K({\bf g})$ acts on
$V_{x}^{{\bf g}}\times {\bf C}^{2}$ through:
$g_{1}(u,v_{1},v_{2})=(u,e^{2\pi i\cdot \frac{2}{3}}v_{1}, e^{2\pi
i\cdot \frac{2}{3}}v_{2})$.

The obstruction bundle $E_{({\bf g})}$ is a plane bundle. And
$dim_{{\bf C}}X_{({\bf g})}=2$, so the 3-point function (5.25) is
nonzero only if $\eta_{j}\in H^{0}(X_{(g_{1})};{\bf Q}),j=1,2,3$.
In this case
$$<\eta_{1},\eta_{2},\eta_{3}>_{orb}=\eta_{1}\eta_{2}\eta_{3}\int_{X_{(\bf{g})}}^{orb}
e_{A}(E_{(\bf{g})}) \eqno{(5.26)}$$ From the section 5.5., the
obstruction bundle $E_{(\bf{g})}$ is the whitney sum of two
orbifold line bundles, let $E_{({\bf g})}=E_{1}\oplus E_{2}$.
$E_{1}$ is generated by  $\xi_{1}\otimes \overline{\omega}$ on the
negiborhood $U_{5}\cap X_{({\bf g})}$; and $E_{2}$ is generated by
$\xi_{2}\otimes \overline{\omega}$ on $U_{5}\cap X_{({\bf g})}$.
So from the first part of this section, consider the orbifold
principal $S^{1}$-bundle $P_{l}$ of  $E_{l}(l=1,2)$. From the
section 5.5, we can see that ${\bf Z}_{d_{l}}=K({\bf g})={\bf
Z}_{3}$, so let  $E_{l}^{'}=(P_{l}/K({\bf g}))\times_{S^{1}}{\bf
C}$ is the orbifold bundle over $X_{({\bf g})}$. $\pi_{K({\bf
g})}:P_{l}\longrightarrow P_{l}/K({\bf g})$ is the projective map.
Note that on every fibre, $\pi_{K({\bf g})}$ is given by:
$z\longmapsto z^{3}$. The lie algebra of $F=S^{1}$ is ${\bf R}$.
So the induced map on the lie algebra is:
 $(\pi_{K({\bf g})})_{*}: {\bf R}\longrightarrow
{\bf R}, a\longmapsto 3a$. From (5.14),
$$\int_{X_{(\bf{g})}}^{orb}
e_{A}(E_{(\bf{g})})=\frac{1}{9}\int_{X_{(\bf{g})}}^{orb}
\Pi_{l=1}^{2}e_{A_{l}^{'}}(E_{l}^{'}) $$ where $A_{l}$ and
$A_{l}^{'}$ are the connections on the bundles $E_{l}$ and
$E_{l}^{'}$ such that $\pi_{K({\bf g})}^{*}(A_{l}^{'})=A_{l}$. In
this moment the group $K({\bf g})$ acts on the bundle $E_{l}^{'}$
trivially, so $E_{l}^{'}$ induced an orbifold bundle $E_{l}^{''}$
over the reduced orbifold $X_{({\bf g})}^{'}$. From  (5.15),
$$\int_{X_{(\bf{g})}}^{orb}
\Pi_{l=1}^{2}e_{A_{l}^{'}}(E_{l}^{'})
=\frac{1}{3}\int_{X_{(\bf{g})}^{'}}^{orb}
\Pi_{l=1}^{2}e_{A_{l}^{''}}(E_{l}^{''}) $$ where $A_{l}^{''}$ is
the connection on the bundle $E_{l}^{''}$ induced from the bundle
$E_{l}^{'}$. Thus

$$\int_{X_{(\bf{g})}}^{orb}
e_{A}(E_{(\bf{g})})=\frac{1}{27}\int_{X_{(\bf{g})}^{'}}^{orb}
\Pi_{l=1}^{2}e_{A_{l}^{''}}(E_{l}^{''})  \eqno{(5.27)}$$

We now calculate the integration $\int_{X_{(\bf{g})}^{'}}^{orb}
\Pi_{l=1}^{2}e_{A_{l}^{''}}(E_{l}^{''})$. The unformizing system
of the bundle $E_{l}^{''}$ over the  reduced orbifold $X_{({\bf
g})}^{'}$ can be described as:

$\bullet$  For $\forall x\in X_{({\bf g})}^{'}$, then $C({\bf
g})/K({\bf g})=1$ is the trivial group, and the action is trivial.

Now we use the localization technique to compute the integration
(5.27). We know that $X_{({\bf g})}^{'}=\overline{O}_{\tau}$ is a
toric variety, $\tau=(v_{0},v_{1},v_{2})$. The three fixed points
by the $T$-action on  $X_{({\bf g})}^{'}$ are $p_{3},p_{4},p_{5}$.
We already computed the $T$-equivariant Euler class of the point
$p_{5}$ at the normal bundle is:
$$e_{T}\left(\nu_{p_{5}}\right)=\lambda_{1}\lambda_{2}$$
Similarly, the $T$-equivariant Euler classes of the points $p_{4}$
and $p_{3}$ at the normal bundles are:
$$e_{T}\left(\nu_{p_{4}}\right)=(\lambda_{1}-\lambda_{2})(-\lambda_{2}),
e_{T}\left(\nu_{p_{3}}\right)=(-\lambda_{1}+\lambda_{2})(-\lambda_{1}).$$

Orbifold line bundle $E_{l}^{''},(l=1,2)$ is trivialized by
$\frac{\partial}{\partial x_{l}}\otimes \overline{\omega},(l=1,2)$
on $U_{5}\cap X_{({\bf g})}^{'}$, $\frac{\partial}{\partial
y_{l}}\otimes \overline{\omega},(l=1,2)$ on $U_{4}\cap X_{({\bf
g})}^{'}$, and $\frac{\partial}{\partial w_{l}}\otimes
\overline{\omega},(l=1,2)$ on $U_{3}\cap X_{({\bf g})}^{'}$. On
the neighborhood  $V_{4}$, we have:
 \begin{eqnarray}
\xi_{1}=\frac{\partial}{\partial x_{1}}=\sum \frac{\partial
y_{j}}{\partial x_{1}}\frac{\partial}{\partial y_{j}}
&=&\frac{\partial}{\partial
x_{1}}(\frac{x_{1}}{x_{4}^{2/3}})\frac{\partial}{\partial
y_{1}} \nonumber \\
&=&y_{5}^{2/3}\frac{\partial}{\partial y_{1}} \nonumber
\end{eqnarray}
Similarly, we can compute $\xi_{2}=\frac{\partial}{\partial
x_{2}}=y_{5}^{2/3}\frac{\partial}{\partial y_{2}}$.

On the neighborhood $V_{3}$, we have:
\begin{eqnarray}
\xi_{1}=\frac{\partial}{\partial x_{1}}=\sum \frac{\partial
w_{j}}{\partial x_{1}}\frac{\partial}{\partial w_{j}}
&=&\frac{\partial}{\partial
x_{1}}(\frac{x_{1}}{x_{3}^{2/3}})\frac{\partial}{\partial
w_{1}}\nonumber \\
&=&w_{5}^{2/3}\frac{\partial}{\partial w_{1}} \nonumber
\end{eqnarray}
Similarly, we have $\xi_{2}=\frac{\partial}{\partial
x_{2}}=w_{5}^{2/3}\frac{\partial}{\partial w_{2}}$. So the
transition function of the orbifold line bundle
$E_{l}^{''},(l=1,2)$
 $$h_{45}:(U_{5}\cap X_{({\bf g})}^{'})\times {\bf C}\supset
 (U_{4}\cap U_{5}\cap X_{({\bf g})}^{'})\times {\bf C}\longrightarrow
 (U_{4}\cap U_{5}\cap X_{({\bf g})}^{'})\times {\bf C}\subset
(U_{4}\cap X_{({\bf g})}^{'})\times {\bf C}$$ is:
 $$h_{45}(x,c)=\left(x,(y_{5}^{2/3})^{3}(x)\cdot c\right) $$
The transition function
 $$h_{35}:(U_{5}\cap X_{({\bf g})}^{'})\times {\bf C}\supset
 (U_{5}\cap U_{3}\cap X_{({\bf g})}^{'})\times {\bf C}\longrightarrow
 (U_{5}\cap U_{3}\cap X_{({\bf g})}^{'})\times {\bf C}\subset
(U_{3}\cap X_{({\bf g})}^{'})\times {\bf C}$$ is:
 $$h_{35}(x,c)=\left(x,(w_{5}^{2/3})^{3}(x)\cdot c\right)$$
So we can define the action of $T$ on the bundle $E_{l}^{''}$ as
$$
\begin{array}{ll} (1)&  t(x,c)=(tx,c)=(tx,\chi^{0}(t)c), t\in T,
(x,c)\in (U_{5}\cap X_{({\bf g})}^{'})\times {\bf C};\\
(2)&  t(x,c)=(tx,y_{5}^{2}(t)c)=(tx,\chi^{-2\rho_{2}}(t)c), t\in
T,
(x,c)\in (U_{4}\cap X_{({\bf g})}^{'})\times {\bf C};\\
(3)&  t(x,c)=(tx,w_{5}^{2}(t)c)=(tx,\chi^{-2\rho_{1}}(t)c), t\in
T, (x,c)\in (U_{3}\cap X_{({\bf g})}^{'})\times {\bf C}.
\end{array}
$$
Then the action of $T$ on $E_{l}^{''}$ at the fixed points
$p_{5},p_{4},p_{3}$  has weights $0,(-2\lambda_{2}),
(-2\lambda_{1})$ respectively.

So from (5.21), we have:

\begin{eqnarray}
\int_{X_{(\bf{g})}^{'}}^{orb}
\Pi_{t=1}^{2}e_{A_{t}^{''}}(E_{t}^{''})&=&\frac{0}{\lambda_{1}\lambda_{2}}
+\frac{(-2\lambda_{2})^{2}}
{(\lambda_{1}-\lambda_{2})(-\lambda_{2})} +
\frac{(-2\lambda_{1})^{2}}
{(-\lambda_{1}+\lambda_{2})(-\lambda_{1})} \nonumber \\
&=&4 \nonumber
\end{eqnarray}
From (5.27),
$$\int_{X_{({\bf g})}}^{orb}e_{A}(E_{({\bf
g})})=\frac{4}{27}.$$ And by (5.26),
$$<\eta_{1},\eta_{2},\eta_{3}>_{orb}=\frac{4}{27}\eta_{1}\eta_{2}\eta_{3}.$$

\subsection*{Acknowledgments}
We thank Professor Banghe Li, Yongbin Ruan ,  jianzhong Pan, David
Cox, and Mainak Poddar for very helpful encouragements. We
especially thank Professor
 Yongbin Ruan for explaining to me to consider this interesting
problem.

\subsection*{}

\bibliographystyle{amsplain}
\bibliography{xibi}

\end{document}